\documentclass{elsarticle}
\usepackage{amssymb, amsmath, amsthm}
\usepackage{graphicx}
\usepackage{epic}
\usepackage{psfrag}

\newtheorem{theorem}{Theorem}[section]
\newtheorem{lemma}[theorem]{Lemma}
\newtheorem{proposition}[theorem]{Proposition}
\newtheorem{corollary}[theorem]{Corollary}

\theoremstyle{definition}

\newtheorem{examples}[theorem]{Examples}

\theoremstyle{remark}
\newtheorem{remark}[theorem]{Remark}

\numberwithin{equation}{section}    

\newcommand\NN{\mathbb{{N}}}
\newcommand\RR{\mathbb{{R}}}
\newcommand\CC{\mathbb{{C}}}
\newcommand\ZZ{\mathbb{{Z}}}

\newcommand\rA{{\rm A}}
\newcommand\ra{{\rm a}}
\newcommand\rB{{\rm B}}
\newcommand\rb{{\rm b}}

\newcommand\rc{{\rm c}}

\newcommand\rd{{\rm d}}
\newcommand\rE{{\rm E}}

\newcommand\rI{{\rm I}}

\newcommand\rT{{\rm T}}

\newcommand\rZ{{\rm Z}}

\newcommand\be{{\boldsymbol e}}

\newcommand\bj{{\bf j}}

\newcommand\bx{{\boldsymbol x}}

\newcommand\bz{{\boldsymbol z}}

\newcommand\beps{{\boldsymbol\varepsilon}}

\newcommand\bI{{\boldsymbol I}}

\newcommand\bP{{\boldsymbol P}}

\newcommand\bX{{\boldsymbol X}}

\newcommand\balpha{{\boldsymbol\alpha}}
\newcommand\bbeta{{\boldsymbol\beta}}

\newcommand\btheta{{\boldsymbol\theta}}

\newcommand\bxi{{\boldsymbol\xi}}

\newcommand\bTheta{{\boldsymbol\Theta}}

\newcommand\bnull{{\boldsymbol0}}
\newcommand\bone{{\boldsymbol1}}

\newcommand\cI{{\mathcal I}}
\newcommand\cJ{{\mathcal J}}

\newcommand\cS{{\mathcal S}}
\newcommand\cU{{\mathcal U}}

\newcommand\nix{\hbox{\phantom{X}}}
\newcommand{\mfrac}[2]%
{\raisebox{0.5pt}{\footnotesize$\dfrac{#1}{#2}$}}
\newcommand{\mbinom}[2]%
{\raisebox{0.5pt}{\footnotesize$\dbinom{#1}{#2}$}}
\def\smmat\{#1&#2\cr#3&#4\}%
   {\!\!\pmatrix{\!#1\!&\!#2\!\cr\!#3\!&\!#4\!}\!\!}
\newcommand\scrm{{\raise0.5pt\hbox{-}}}

\begin{document}

\title{Scalar Multivariate Subdivision Schemes\\and Box Splines}

\author[ch]{Maria Charina\corref{cor1}}
\ead{maria.charina@uni-dortmund.de}
\author[cc]{Costanza Conti}
\ead{costanza.conti@unifi.it}
\author[hoh]{Kurt Jetter}
\ead{Kurt.Jetter@uni-hohenheim.de}
\author[hoh]{Georg Zimmermann}
\ead{Georg.Zimmermann@uni-hohenheim.de}

\cortext[cor1]{Corresponding author}
\address[ch]{Fakult\"at f\"ur Mathematik, TU Dortmund, D--44221 Dortmund, Germany}
\address[cc]{Dipartimento di Energetica, Universit\`a di Firenze, Via C. Lombroso 6/17, I--50134 Firenze, Italy}
\address[hoh]{Institut f\"ur Angewandte Mathematik und Statistik, Universit\"at Hohenheim, D--70593 Stuttgart, Germany}

\begin{abstract}

We study scalar $d$-variate subdivision schemes, with dilation matrix $2I$,
satisfying the sum rules of order $k$. Using the
results of M\"oller and Sauer, stated for general
expanding dilation matrices, we characterize the structure of the mask
symbols of such schemes by showing that they must be linear
combinations of shifted box spline generators of some quotient polynomial
ideal. The directions of the corresponding box splines are
columns of certain unimodular matrices. The quotient
ideal is determined by the given
order of the sum rules or, equivalently, by the order of the
zero conditions.

The results presented in this paper open a way to a systematic study
of subdivision schemes, since  box spline subdivisions turn out to
be the building blocks of any reasonable multivariate subdivision scheme.

As in the univariate case, the characterization we give  is the
proper way of matching the smoothness of the box spline building
blocks with the order of polynomial reproduction of the
corresponding subdivision scheme. However, due to the interaction
of the building blocks, convergence and smoothness properties may
change, if several convergent schemes are combined.

The results are illustrated with several examples.
\end{abstract}


\begin{keyword} 
      Subdivision schemes, box splines, ideals
\end{keyword}
\maketitle
\section*{Introduction}  

Subdivision schemes are efficient iterative procedures for
generating finer and finer grids of points in $\RR^d$ and are used
to design smooth curves or surfaces. Starting with some initial
grid of points, a binary scalar subdivision scheme computes the coordinates
of the finer grid points $\rd^{(r+1)}$, inheriting the
topology of the coarser ones, via local averages
\begin{equation}\label{eq:subrecursion}
  \rd^{(r+1)}=\cS_\ra \rd^{(r)}=\sum_{\bbeta\in\ZZ^d}
  \rd^{(r)}_\bbeta\;\ra_{\balpha-2\bbeta}, \quad r \ge 0\;.
 \notag
\end{equation}
These averaging rules depend on the coefficients of the
corresponding subdivision mask $\ra=\big(\ra_\balpha \big)_{\balpha\in \ZZ^d}$
which we assume to be some finitely supported sequence of real
numbers. If the mask is chosen appropriately, the grids with the
vertices $\rd^{(r)}$, as $r$ goes to infinity, can be interpreted
as approximations of the values of a limiting curve or surface.
The locality of the method and its algorithmic simplicity ensure
that the subdivision recursion is fast, efficient, and easy to
implement. These features and the connection between subdivision
and other multiresolution methods have led to an increasing
popularity of subdivision in computer graphics, in computer aided
geometrical design, and in wavelet and frame constructions. For
more details on subdivision we refer the interested reader to the
pioneering work in \cite{CDM91, D92} or to the more recent survey
\cite{DynLevin02} and the references therein. For basic details on
wavelet and frame constructions see, {\it e.g.},
\cite{Chr,Daub,StrNgu}.

\medskip
In general, the limit of a convergent subdivision scheme is not known
analytically. Nevertheless, various analytical properties of the limit
can be read off the mask symbol
\begin{equation*}
  \ra(\bz)=\sum_{\balpha \in \ZZ^d} \ra_\balpha \bz^\balpha\,,
  \qquad \bz \in \big(\CC \setminus \{0\}\big)^d\;,
\end{equation*}
which---due to the finite support of the mask---is a
$d$-variate Laurent polynomial. In the univariate setting, the
properties of convergent subdivision schemes are well-understood.
The symbol $\ra(z)$ of any convergent univariate subdivision
scheme possesses the \underline{Factorization Property}: it can be
written as
\begin{equation*}\label{eq:b-splines}
  \ra(z)=2\cdot \sigma(z)\cdot \Bigl(\mfrac{1+z}{2}\Bigr)^k\;,
\quad\text{with}\quad
  \sigma(1)=1,
\end{equation*}
for some integer $k \ge 1$, with the Laurent polynomial
$\sigma(z)$  satisfying some additional properties to ensure the
convergence and smoothness of the subdivision limit, see \cite{D1,
Plonka}. Since the coefficients of $\sigma$ sum to $1$, this
Factorization Property tells us that convergence of the
subdivision scheme implies that its mask can be written as an
affine combination of shifted versions of B-spline symbols
$B_k(z)=(1+z)^k / 2^{k-1}$. Thus, any reasonable univariate
subdivision scheme uses B-spline symbols as its building blocks, a
fact which is very crucial for designing new efficient subdivision
processes, see {\it e.g.}, \cite{ContiRomani, DD, DongShen,
DynHormSS}.

\medskip
In this paper we are going to present results which aim at
replacing the Factorization Property of univariate schemes by a
\underline{Decomposition Property} of multivariate schemes, see
Theorems~A--C below. Theorem~C  deals with the case $d=2$ and is
crucial for studying the properties of subdivision schemes. It
tells us that the mask symbol of any reasonable bivariate
subdivision scheme can be decomposed as
\begin{equation*}\label{eq:box-spline}
  \ra(z_1,z_2)=4 \cdot \sum_{B^{\#}_{\alpha,\beta,\gamma}\in \rI_k}
  \lambda_{\alpha,\beta,\gamma} \cdot
  \sigma_{\alpha,\beta,\gamma}(z_1,z_2) \cdot
  B_{\alpha,\beta,\gamma}^\# (z_1,z_2)\;,
\end{equation*}
where $\displaystyle{\sum\lambda_{\alpha,\beta,\gamma}=1}$\,, and
the symbols
 $\sigma_{\alpha,\beta,\gamma}(z_1,z_2)$ are Laurent polynomials
normalized by $\sigma_{\alpha,\beta,\gamma}(1,1)=~1$. Here,
\begin{equation*}
  B^{\#}_{\alpha,\beta,\gamma}(z_1,z_2) =
  \Bigl(\mfrac{1+z_1}{2}\Bigr)^\alpha
  \Bigl(\mfrac{1+z_2}{2}\Bigr)^\beta
  \Bigl(\mfrac{1+z_1\,z_2}{2}\Bigr)^\gamma, \quad
  \alpha,\beta,\gamma \in \NN_0\;,
\end{equation*}
are the normalized mask symbols of three-directional box splines,
$k$ refers to the order of polynomial reproduction of the
subdivision operator $S_\ra$---in the wavelet literature also
denoted as order of accuracy of the mask---and the list $\rI_k$ is
as in Theorem~\ref{thm:twoDJkgen}. The advantage of such a
decomposition result is obvious, not only for the purpose of
classifying the zoo of subdivision schemes, but also as a starting
point for constructing new schemes or for enhancing the properties
of the schemes according to specific requirements.

\medskip
We believe that the use of the building blocks
$B^{\#}_{\alpha,\beta,\gamma}$ in the above decomposition  is the
appropriate generalization to the bivariate case of the normalized
univariate B-spline symbols, since
$B^{\#}_{\alpha,\beta,\gamma}$ correspond to a class of
well-known bivariate spline functions whose order of global
smoothness and whose order of polynomial reproduction match in the
same way as in the univariate case, see the notes at the end of
Section~\ref{ssec:prop}. We consider, therefore, this paper to be
also an interesting contribution to the box spline literature.

\medskip
The methods employed in this paper are of algebraic nature, as
well as other results on subdivision dealing with polynomial
reproduction, see \cite[Chapter~6]{CDM91}. It is well-known that
the necessary conditions for the convergence of a subdivision
scheme are equivalent to the fact that  the mask symbol $\ra(\bz)$
is a properly normalized element of a certain ring of Laurent
polynomials: the symbol must belong to the ideal $\cI$, or to its
power $\cI^k$, depending on the order $k$ of the polynomial
reproduction of the scheme. For definitions of $\cI^k$, $k \in
\NN$, see  \eqref{eq:ideal} and \eqref{eq:sumruleone}, or see the
zero conditions of order $k$ in \eqref{eq:Zkcond}. This algebraic
property of $\ra(\bz)$ has been studied in detail in the papers by
Sauer \cite{Sau02a, Sauer} and by M\"oller and Sauer
\cite{MoeSau}, which also motivated us to consider these ideals
once again. Note that the Factorization Property says that
$\ra(\bz)$ belongs to the principal ideal generated by some
B-spline symbol, while the Decomposition Property tells us that
$\ra(\bz)$ is in the ideal generated by a particular set of box
spline symbols. M\"oller and Sauer have given other sets of
generators for these ideals, putting the emphasis on the algebraic
properties of these generators. We consider it one of our main
achievements that we were able to relate the generators of the
ideals $\cI^k$ to a well-known class of spline functions, see
\cite{BH2}. Box spline subdivision has been studied thoroughly in
the literature, see \cite{BHR,Chui}, where such box spline schemes
also appear under a different name. In the three-directional case,
for example, their mask symbols are given by $4 \cdot
B^{\#}_{\alpha,\beta,\gamma}(z_1,z_2)$.

\medskip
Starting with Section 2, we restrict ourselves to working with the
ideals of the ring of $d$-variate polynomials instead of the ring
of Laurent polynomials. Note that any Laurent polynomial can be
shifted to produce a polynomial. It is done by multiplying the
Laurent polynomial with a factor $\bz^\balpha$, {\it i.e.}, a unit
in the ring of Laurent polynomials.  This results in a shift of
the support of the mask. Such a shift does neither affect the
convergence and regularity of the  scheme, nor does it change the
zero conditions of any order. We, therefore, assume that the mask
is supported in $\NN_0^d$.

\medskip
The paper is organized as follows: in Section 1 we introduce
some notation and background on scalar   subdivision
schemes with dilation matrix $2I$. In Section 2 we show how the  zero
conditions of order $k=1$ on
the mask symbol $\ra(\bz)$ determine the structure of
the symbol of any convergent binary subdivision scheme. This result is given
in Theorem~A
and is a crucial step toward understanding the case {$k>1$}. As in the
multivariate case the Factorization Property is replaced by
the Decomposition Property, it is important to
reduce the number of the
terms of this decomposition or, equivalently, to minimize the number of
the generators for $\cI$, see Theorem~\ref{thm:JBoxgen},
Theorem~\ref{thm:JBoxgenmod}.
Section 3 contains our main results, Theorems~B and~C, dealing with powers
of $\cI$. The fact that the
ideal $\cI^k$ is generated
by the appropriate
products of the elements of $\cI$ makes the proof of Theorem~B
straightforward, if one
is not interested in reducing the number of terms of the corresponding
decomposition.
We address the
latter issue only in the case $d=2$, as this case is of special interest in
subdivision. The
corresponding result is stated in Theorem~C.
In Section~4 we illustrate the result of Theorem~C  with several bivariate
examples.

\medskip
This paper is an extended version of the technical report
\cite{ccjzrep} which can be obtained from the first author.

\section{Background on subdivision and notation} 
A scalar $d$-variate subdivision scheme is given by a scalar
$\ZZ^d$-indexed sequence $\ra = (\ra_\balpha)_{\balpha\in\ZZ^d}$,
the so-called mask, defining the subdivision operator $\cS_\ra$ on
data sequences $\rd=(\rd_\balpha)_{\balpha\in\ZZ^d} \in
\ell(\ZZ^d)$ as follows:
\begin{equation}\label{eq:subdop}
  \big(\cS_\ra \rd\big)_\balpha = \sum_{\bbeta\in\ZZ^d}
  \rd_\bbeta\;\ra_{\balpha-2\bbeta}\;,\quad \balpha\in\ZZ^d\;.
\end{equation}
We assume that the mask is finite, {\it i.e.}, only finitely many
coefficients $\ra_{\balpha}$ are non-zero.

In our study we use the following symbol notation. For a
finitely supported sequence $\rc=(\rc_\balpha)_{\balpha\in\ZZ^d}$,
its symbol is given by the Laurent polynomial
\begin{equation*}\label{eq:laurpol}
  \rc(\bz) = \sum_{\balpha\in\ZZ^d} \rc_\balpha\;\bz^\balpha
\end{equation*}
with $\bz = (z_1,\dots,z_d) \in\left( \CC\setminus\{0\}\right)^d$
and, in the multi-index notation,
\begin{equation*}
  \bz^\balpha = z_1^{\alpha_1}z_2^{\alpha_2}\cdot \ldots \cdot z_d^{\alpha_d}\;,
  \quad \text{for } \balpha = (\alpha_1,\dots,\alpha_d) \in \ZZ^d\;.
\end{equation*}
In the symbol notation, the subdivision step in \eqref{eq:subdop}
is described by the identity
\begin{equation}\label{eq:subdopsymb}
  \big(\cS_\ra \rd\big)(\bz) = \rd(\bz^2)\;\ra(\bz)\;, \quad
\bz^2 = (z_1^2, z_2^2, \dots, z_d^2).
\end{equation}
The first factor on the right-hand side of  \eqref{eq:subdopsymb} refers to
an upsampled
version of the data $\rd$.

Equation \eqref{eq:subdopsymb} can also be written
using convolution with the so-called submasks of $\ra$. Let
\begin{equation}\label{eq:extrpts}
  \rE = \{0\,{,}\,1\}^d
\end{equation}
be the set of representatives of $\ZZ^d/2\ZZ^d$, given by the
vertices of the unit cube $[0,1]^d$, containing
\begin{equation}\label{eq:eone}
  \bnull = (0,0,\dots,0) \quad\text{and}\quad \bone = (1,1,\dots,1)\,.
\notag
\end{equation}
Then, the $2^d$ submasks $\ra_\be$ and their symbols $\ra_\be(\bz)$
are defined by
\begin{equation} \label{eq:subsymbols}
  \ra_\be = (\ra_{\be+2\balpha})_{\balpha\in\ZZ^d} \quad \text{and}
  \quad \ra_\be(\bz) = \sum_{\balpha\in\ZZ^d}
  \ra_{\be+2\balpha}\;\bz^\balpha\;,\quad \be\in\rE\;. \notag
\end{equation}
The standard decomposition
\begin{equation} \label{eq:decomp}
  \ra(\bz) =
  \sum_{\be\in\rE} \bz^{\be}\;\ra_{\be}(\bz^2) \notag
\end{equation}
yields the equivalent form of the identity
\eqref{eq:subdopsymb}
\begin{equation} \label{eq:subrules}
  \big(\cS_\ra \rd\big)(\bz) = \sum_{\be\in\rE} \bz^\be \,
  \rd(\bz^2)\,\ra_\be(\bz^2)\;. \notag
\end{equation}
This shows that a subdivision step is the result of convolving the input
data $\rd$ with each submask $\ra_{\be}$, which is a process of low pass
filtering, followed by an interleaving process, {\it i.e.}, upsampling and
multiplication
by $\bz^\be$, to produce the output data $\cS_\ra \rd$.

\medskip
We say that the subdivision scheme $\cS_a$ is convergent, if for
any starting sequence $\rd \in \ell_\infty(\ZZ^d)$, there exists a
uniformly continuous function $f_\rd$ such that
\begin{equation*}
  \lim_{r \rightarrow \infty} \sup_{\balpha \in \ZZ^d} \left|
  (\cS_a^r \rd)_\balpha - f_\rd(2^{-r} \balpha)\right|=0 \;,
\end{equation*}
and $f_\rd \neq 0$ for some initial data $\rd$. This is the notion
of $C$-convergence, also referred to as uniform convergence, see
\cite{CDM91}, where $L_p$-convergence for $1 \le p < \infty$ is
also discussed. The necessary condition for these types of
convergence is now known to be the so-called sum rule of order~$1$
referring to the submasks
\begin{equation}\label{eq:neccond}
  \ra_\be(\bone) = \sum_{\balpha\in\ZZ^d} \ra_{\be+2\balpha} = 1\;,
  \quad \be\in\rE\;,
\end{equation}
see \cite[Proposition 2.1]{CDM91} and \cite[Theorem 3.1]{HJ98}.

Another notion of convergence is used in the literature on
multiresolution methods devoted to wavelet and frame
constructions, for details see \cite[Chapter~13]{Chr}, 
\cite[Chapter~7]{Daub}, or \cite[Chapter~6]{StrNgu}. There, the
convergence is characterized by the properties of the infinite
product
$$
\prod_{j=0}^\infty \ra^\#(\bz^{2^j}) \;, \quad
\ra^\#(\bz) = \frac{1}{2^d}\; \ra(\bz),
$$
if $\bz$ is restricted to the $d$-dimensional torus. To state
some of the properties we switch to the real variables
$\bxi=(\xi_1,\dots,\xi_d)$ via the transformation
$z_j=e^{-i\pi\xi_j}$, $j=1,\dots,d$. The set \eqref{eq:extrpts}
then transforms into the set
\begin{equation}\label{eq:extrptstwo}
  \rZ=\rZ_\rE = \{\beps =  e^{-i\pi\be} \;:\; \be\in \rE\} =
  \{{-}1\,{,}\,{+}1\}^d
\end{equation}
of the vertices of the cube $[-1,+1]^d$, and the necessary condition
\eqref{eq:neccond} takes the equivalent form
\begin{equation} \label{eq:sumruleone}
 \ra(\bone) = 2^d \quad \text{and} \quad \ra(\beps) = 0 \quad
  \text{for }\beps \in \rZ'=\rZ\setminus \{\bone\}\;.
\end{equation}
For this reason we call $\rZ'$ the zero set, and the conditions in
\eqref{eq:sumruleone} the \emph{zero condition of order one
(Condition~Z$_1$)}. In the literature, both the conditions in
\eqref{eq:neccond} and their equivalent form
 in \eqref{eq:sumruleone} are   called the sum rules of order one.

\medskip

More generally, we  also use the higher order sum rules following
the notation introduced and discussed in the survey paper
\cite{JePlo}, see also the references therein: The mask symbol
$\ra(\bz)$ is said to satisfy the \emph{zero condition of order
$k$ (Condition~Z$_k$)}, if
\begin{equation} \label{eq:Zkcond}
\ra(\bone)=2^d \quad \hbox{and} \quad
\bigl(D^\bj \ra\bigr)(\beps)=0 \quad
\hbox{for} \quad \beps\in \rZ' = \rZ\setminus\{\bone\}\quad \hbox{and} \quad
|\bj|<k\;.
\end{equation}

\section{Zero condition and the associated ideal} 
\label{sec:zeroconideal}

In this section we show that condition Z$_1$
 fully determines the structure of
 the symbol
$\ra(\bz)$ of any convergent scalar subdivision scheme, see Theorem~A.
The statement of this theorem is therefore the first step toward   the desired
multivariate
generalization of the Factorization Property of univariate schemes.
Theorem~A is also  crucial for understanding the structure of the symbol
 $\ra(\bz)$ satisfying condition~Z$_k$ for $k>1$. It  also
shows  that
the Factorization Property is replaced in the multivariate case by the
decomposition \eqref{eq:combA}.

The generators $q_\bTheta$ in \eqref{eq:combA} are the symbols of
subdivision schemes whose limit functions $1_\Theta$ are the
characteristic functions of the parallelepiped spanned by the
column vectors of certain matrices $\bTheta$. These functions are
box splines of degree zero. In order to control the size of the
mask, it is, thus, of importance to choose these box spline
symbols appropriately, and to work in \eqref{eq:combA} with as few
summands as possible. Theorem~\ref{thm:JBoxgen} shows how to meet
these requirements.

A modification of the generators $q_\bTheta$ is given in
Theorem~\ref{thm:JBoxgenmod}
and Theorem~$\widetilde A$, at the end of this section. The modified
generators $\widetilde{q}_\bTheta$ have a useful  algebraic property:
they factor
into $d$ linear polynomials.

\medskip
Since all convergent subdivision schemes satisfy condition Z$_1$,
we start with a characterization of the polynomial ideal
\begin{equation} \label{eq:ideal}
  \cI = \bigl\{ p\in\Pi^d : p(\beps)=0
  \text{ for } \beps\in \rZ' \bigr\}\;,
\end{equation}
and later, in Section~3, of its powers
\begin{equation*}
  \cI^k = \bigl\{ p\in\Pi^d : \bigl(D^\bj p\bigr)(\beps)=0
  \text{ for } \beps\in \rZ'\,,\; |\bj|<k
  \bigr\}
\end{equation*}
for $k > 1$. The main result of this section, whose proof is a consequence of
Theorem~\ref{thm:JBoxgen}, states the following:

\medskip
\noindent
{\bf Theorem A.} {\it
The mask symbol of any convergent $d$-variate subdivision scheme
$\cS_\ra$ can be written in the form
\begin{equation} \label{eq:combA}
  \ra(\bz) = \sum_{\bTheta} \lambda_\bTheta\; \sigma_\bTheta(\bz)\;
  2^d\,q_\bTheta(\bz)\;.
\end{equation}
The sum runs over
all unimodular $d{\times}d$-submatrices $\bTheta$ of $\bigl(\bX_d^{(1)}
\setlength{\unitlength}{1mm}\begin{picture}(2,3)
\put(0,0){\dottedline{0.5}(1,-1)(1,3)}
\end{picture} \bX_d^{(2)}\bigr)$ from
\eqref{eq:directions}. The polynomials $q_\bTheta(\bz)$ are
defined in \eqref{eq:generators}, $\sigma_\bTheta(\bz)$ are
Laurent polynomials satisfying $\sigma_\bTheta(\bone)=1$\,, and
$\lambda_\bTheta$ are real numbers subject to $\sum_{\bTheta}
\lambda_\bTheta = 1$\,.}

\begin{remark}
Since the $2^d\,q_\bTheta(\bz)$ are the mask symbols of certain box
splines of degree zero, the mask $\ra$ is an affine combination of
masks each of which originates from such a box spline convolved with
some (smoothing) factor. Thus, Theorem~A explains why such affine
combinations were successfully studied before, see {\it e.g.}, \cite{Conti,
CGPS07, ContiRomani} for examples of bivariate and univariate schemes.
\end{remark}

From the point of view of  algebraic geometry, Theorem~A tells us that
the system of box spline symbols $q_\bTheta$ in the representation
\eqref{eq:combA} generates the ideal $\cI$. This connects our work to
the papers \cite{Sau02a, Sauer} where the author has studied these
ideals in detail and has characterized the ideal $\cI$ in
\cite[Proposition 4.1]{Sau02a} using the generators
\begin{equation}\label{eq:genideal}
  \cI = {<}\, z_1^2-1, z_2^2-1, \ldots,z_d^2-1,
  (z_1+1)(z_2+1)\ldots(z_d+1) \,{>}\,.
\end{equation}
The characterization \eqref{eq:genideal}
 is  a special case of his more general results,
see also \cite[Example 4]{MoeSau}. It also follows from the
fact that the polynomial ideal
\begin{equation*}
  \cJ := \bigl\{ p\in\Pi^d : p(\beps)=0
    \text{ for } \beps\in \rZ \bigr\}
\end{equation*}
of polynomials vanishing on $Z$ is generated by the polynomials
$$
  z_1^2-1, \ z_2^2-1, \ \ldots, \ z_d^2-1\,.
$$
This result is stated in \cite[Lemma 2.3]{CDM91} with an elementary and
constructive proof. Consequently, in order to determine a set of
generators for the quotient ideal $\cI$, it suffices to add to the generators of $\cJ$
the polynomial 
$$
 (z_1+1)(z_2+1)\ldots(z_d+1)
$$ 
that vanishes on $\rZ'$\,, but does not vanish at $\bone$\,.

\begin{figure}
\psfrag{ze}[cc][cc]{$z_1$}
\psfrag{zz}[cc][cc]{$z_2$}
\psfrag{zep1io}[cc][cc]{$1{+}z_1=0$}
\psfrag{zzp1io}[cc][cc]{$1{+}z_2=0$}
\psfrag{zepzzio}[cc][cc]{$z_1{+}z_2=0$}
\psfrag{zemzzie}[cc][cc]{$1{+}z_1\,z_2=0$}
\psfrag{0}[cc][cc]{$0$}
\psfrag{e}[cc][cc]{$1$}
\psfrag{me}[cc][cc]{$-1$}
$$\includegraphics{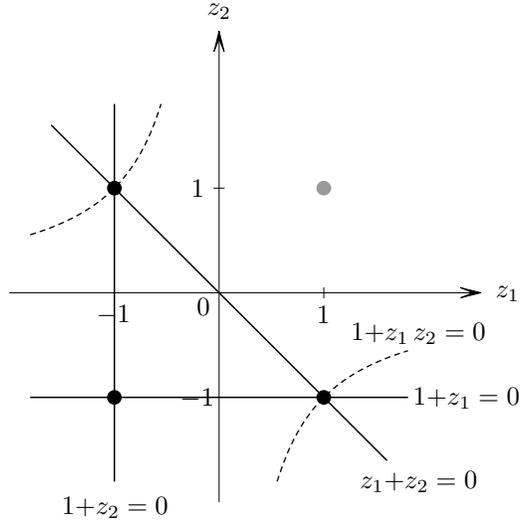}$$
\caption{The ideal $\cI$ from the algebraic geometer's point of view.}
\label{fig:alggeom}
\end{figure}

\medskip
In the two-dimensional case, guided by ideas from algebraic geometry,
it seems natural to use the straight lines through the three points of
$\rZ'$, as shown in Figure~\ref{fig:alggeom}, in order to find other
generators for the ideal $\cI$\,. This yields the following result.

\hspace{3cm}

\begin{proposition}\label{prp:Jd2Boxmod}
For $d=2$\,, we have
\begin{equation}\label{eq:Jd2Boxmod}
  \cI = {<}\, \mfrac{1+z_1}{2}\,\mfrac{1+z_2}{2}   \,,
              \mfrac{1+z_1}{2}\,\mfrac{z_1+z_2}{2} \,,
              \mfrac{1+z_2}{2}\,\mfrac{z_1+z_2}{2} \,{>} \,.
\end{equation}
\end{proposition}

\begin{proof}
By~\eqref{eq:genideal}, we have for $d=2$
\begin{equation*}
  \cI = {<}\, z_1^2-1, \ z_2^2-1,\ (z_1+1)(z_2+1) \,{>} \,.
\end{equation*}
Let us denote
\begin{equation*}
  \widetilde\cI = {<}\, \mfrac{1+z_1}{2}\,\mfrac{1+z_2}{2}   \,,
                        \mfrac{1+z_1}{2}\,\mfrac{z_1+z_2}{2} \,,
                        \mfrac{1+z_2}{2}\,\mfrac{z_1+z_2}{2} \,{>} \,,
\end{equation*}
then it suffices to show that the generators for $\cI$ are contained in
$\widetilde\cI$ and \textit{vice versa}.

To this end, note first that
\begin{align*}
  z_1^2-1 &= -4\,\mfrac{1+z_1}{2}\,\mfrac{1+z_2}{2}
             +4\,\mfrac{1+z_1}{2}\,\mfrac{z_1+z_2}{2}\\
  z_2^2-1 &= -4\,\mfrac{1+z_1}{2}\,\mfrac{1+z_2}{2}
             +4\,\mfrac{1+z_2}{2}\,\mfrac{z_1+z_2}{2}\\
  (z_1+1)(z_2+1) &= 4\,\mfrac{1+z_1}{2}\,\mfrac{1+z_2}{2} \,,
\end{align*}
and conversely that
\begin{align*}
  \mfrac{1+z_1}{2}\,\mfrac{1+z_2}{2} &= \mfrac{1}{4}\,(z_1+1)(z_2+1)\\
  \mfrac{1+z_1}{2}\,\mfrac{z_1+z_2}{2} &=
    \mfrac{1}{4}\,(z_1^2-1)+\mfrac{1}{4}\,(z_1+1)(z_2+1)\\
  \mfrac{1+z_2}{2}\,\mfrac{z_1+z_2}{2} &=
    \mfrac{1}{4}\,(z_2^2-1)+\mfrac{1}{4}\,(z_1+1)(z_2+1) \,.
\end{align*}
This shows that $\widetilde\cI = \cI$ as claimed.
\end{proof}

\begin{remark}\label{rem:gridrefl}
The three functions in~\eqref{eq:Jd2Boxmod} are remarkably close to the
symbols of the three box splines of degree zero on the three-directional
grid. Indeed, if we reflect the standard three-directional
grid about one of the coordinate axes, \textit{i.\,e.}, use the grid
spanned by the three vectors
$\left(\begin{smallmatrix}1\\[0.5mm]0\end{smallmatrix}\right),
 \left(\begin{smallmatrix}0\\[0.5mm]1\end{smallmatrix}\right)$, and
$\left(\begin{smallmatrix}1\\[0.5mm]-1\end{smallmatrix}\right)$\,,
then the box splines of degree zero have the mask symbols
\begin{equation*}
  4\,\mfrac{1+z_1}{2}\,\mfrac{1+z_2}{2}     \,,\qquad
  4\,\mfrac{1+z_1}{2}\,\mfrac{1+z_1/z_2}{2} \,,\qquad
  4\,\mfrac{1+z_2}{2}\,\mfrac{1+z_1/z_2}{2} \,.
\end{equation*}
Since we can write $z_1+z_2 = z_2\,(1+z_1/z_2)$\,,
these are just the functions in~\eqref{eq:Jd2Boxmod} up to
an appropriate normalization.

It follows that the result of Proposition \ref{prp:Jd2Boxmod}
can as well use
the usual box spline symbols: The reflection of
the three-directional grid about the $z_1$-axis, as mentioned before,
corresponds to the variable transformation $(z_1\,{,}\,z_2)\mapsto
(z_1\,{,}\,1/z_2)$\,, and then~\eqref{eq:Jd2Boxmod} becomes---again
after a proper normalization---
\begin{equation}\label{eq:Jd2Box}
  \cI = {<}\, \mfrac{1+z_1}{2}\,\mfrac{1+z_2}{2}      \,,
              \mfrac{1+z_1}{2}\,\mfrac{1+z_1\,z_2}{2} \,,
              \mfrac{1+z_2}{2}\,\mfrac{1+z_1\,z_2}{2} \,{>} \,.
\end{equation}
From the algebraic geometer's point of view, this amounts to replacing
the straight line $z_1+z_2=0$ in Figure~\ref{fig:alggeom} by the
hyperbola $1+z_1\,z_2=0$\,.
\end{remark}

\medskip
The identity \eqref{eq:Jd2Box} is the statement of
Theorem~\ref{thm:JBoxgen} for the case $d=2$\,. To be able to
state this result in general, we first need to provide some
additional notation. As before, we denote by $\rE =
\{0\,{,}\,1\}^d$ the set of vertices of the $d$-dimensional
hypercube, and let $\rE'=\rE\setminus\{\bnull\}$\,. Collecting the
elements of $\rE'$ in the  matrix $ \bX_d = \Big( \;\;\be\;\; \Big)_{\be\in\rE'}$ yields 
\begin{gather}\label{eq:directions}
  \setlength{\unitlength}{1mm}%
  \begin{picture}(0,0)
  \put(0,0){\dottedline{0.5}(37,-12)(37,12)}
  \put(0,0){\dottedline{0.5}(80,-12)(80,12)}
  \put(0,0){\dottedline{0.5}(95,-12)(95,12)}
  \end{picture}%
  \bX_d =
  \setcounter{MaxMatrixCols}{20} 
  \begin{pmatrix}
  1      & &        & & 0      & 1      & 1      &        & 1      & 0
    &        & 0      & &        & & 1\\[-2mm]
  0      & &        & & \vdots & 1      & 0      &        & 0      & 1
    &        & \vdots & &        & & 1\\[-2mm]
  \vdots & & \cdots & & \vdots & 0      & 1      & \cdots & \vdots & 1
    & \cdots & 0      & & \cdots & & \vdots\\[-2mm]
  \vdots & &        & & 0      & \vdots & 0      &        & 0      & 0
    &        & 1      & &        & & 1\\[-2mm]
  0      & &        & & 1      & 0      & \vdots &        & 1 & \vdots
    &        & 1      & &        & & 1
  \end{pmatrix}\,,\\ \notag
  \text{\footnotesize$\hspace{11mm} \bX_d^{(1)}=\bI_d
  \hspace{24mm} \bX_d^{(2)}
  \hspace{21mm} \cdots
  \hspace{6mm}  \bX_d^{(d)}$}
\end{gather}
where each column of the submatrix $\bX_d^{(k)}$ contains exactly
$k$ entries equal to~$1$\,. Thus $\bX_d^{(1)}$ consists of the standard
unit vectors $\be_k$\,, $k=1,\dots,d$, while $\bX_d^{(2)}$ contains the
vectors $\be_j+\be_k$, $j\not=k$, \textit{etc}.
We treat the columns of $\bX_d$ as
directional vectors from which we build the box splines of degree zero,
\textit{i.\,e.}, the characteristic functions of certain
parallelepipeds. This means that we take any $d$ columns from $\bX_d$
to produce a square submatrix $\bTheta$ of $\bX_d$\,. With each such
$\bTheta$\,, we associate the normalized polynomial
\begin{equation}\label{eq:generators}
  q_\bTheta (\bz) = \prod_{\btheta \in \bTheta}
  \mfrac{1+\bz^\btheta}{2}\;,
\end{equation}
where $\btheta$ runs through the columns of $\bTheta$\,, with
$q_\bTheta (\bone)=1$\,.

\begin{proposition}\label{prp:nonunimodcase}
For any $d{\times}d$-submatrix $\bTheta$ of $\bX_d$ in
\eqref{eq:directions}, we have
\begin{equation*}
  q_\bTheta \in \cI \quad \Longleftrightarrow \quad
  \det \bTheta \equiv 1 \;(\operatorname{mod}\; 2)\;.
\end{equation*}
\end{proposition}

\begin{proof}
By definition, $q_\bTheta \in \cI$ if and only if
\begin{equation*}
  \forall \beps\in \rZ'\,: \qquad
  q_\bTheta(\beps) = \prod_{\btheta\in\bTheta}
  \mfrac{1+\beps^\btheta}{2} = 0 \,.
\end{equation*}
Due to $\beps^\btheta =
e^{-i\pi\be^T\btheta}$\,, $\beps \in \rZ'$,  this is equivalent to
\begin{equation*}
  \forall\be\in \rE'\,:\; \exists \btheta\in\bTheta\,: \qquad
  \be^T\btheta \equiv 1 \;(\operatorname{mod}\; 2)\;.
\end{equation*}
\textit{I.\,e.}, the map
\begin{equation*}
  L\,:\; \ZZ^d \to \ZZ^d\,,\quad \bx^T \mapsto \bx^T\bTheta\,,
\end{equation*}
has the property
\begin{equation} \label{eq:Lprop}
  \forall\be\in \rE'\,: \qquad
  L(\be^T) \not\equiv \bnull^T \;(\operatorname{mod}\; 2)\;.
\end{equation}

As we are only interested in the parity of $\det \bTheta$\,, we may
employ the ring homomorphism $\ZZ\to\ZZ_2=\ZZ/2\,\ZZ$ mapping even
and odd numbers to their coset representatives $\overline{0}$ and
$\overline{1}$\,, respectively. It extends naturally to a matrix ring
homomorphism $\ZZ^{d\times d} \to (\ZZ_2)^{d\times d}$ and thus induces
a linear map $\overline{L} : (\ZZ_2)^d \to (\ZZ_2)^d$ satisfying
$\det \overline{L} = \overline{\det L}$\,.
Since the set $\rE$ is a complete set of coset representatives of the
vector space $(\ZZ_2)^d = \ZZ^d / 2\ZZ^d$\,, property \eqref{eq:Lprop}
is equivalent
to
\begin{equation*}
  \ker (\overline{L}) = \{\overline{\bnull}^T\}\,,
\end{equation*}
\textit{i.\,e.}, $\overline{L}$ is a vector space automorphism and thus
\begin{equation*}
  \det \overline{L} \not= \overline{0} \quad\text{in }\ZZ_2\,.
\end{equation*}
In other words,
\begin{equation*}
  \det \bTheta = \det L \not\equiv 0 \;(\operatorname{mod}\; 2)\;,
\end{equation*}
which yields the claim.
\end{proof}

The next result, which is important for the proof of Theorem~A,
states that, conversely, the ideal $\cI$ is generated by a
subfamily of the elements $q_\bTheta$. It is worth noting that the
condition $\det \bTheta = \pm1$ is equivalent to $\bTheta$ having
an integer inverse and, thus, to the fact that its columns
generate the integer grid $\ZZ^d$\,.

\begin{lemma}\label{lem:Jgensingle}
For $k=1,\dots,d$\,, we have
\begin{equation*}
  z_k^2-1 \;\in\; {<}\,\Bigl\{ q_\bTheta : \bTheta \;
  d{\times}d\text{-submatrix of } \bigl(\bX_d^{(1)}
  \setlength{\unitlength}{1mm}\begin{picture}(2,3)
  \put(0,0){\dottedline{0.5}(1,-1)(1,3)}
  \end{picture} \bX_d^{(2)}\bigr)\,,\;\det \bTheta = \pm1\,,\;
  \be_k \in\bTheta \Bigr\}\,{>} \,.
\end{equation*}
\end{lemma}

\begin{proof}
For $d=1$, the lemma simply claims that
\begin{equation*}
  z^2-1 \;\in \;{<}\, \mfrac{1+z}{2} \,{>} \,,
\end{equation*}
which follows from $z^2-1 = 2\,(z-1)\,\mfrac{1+z}{2}$\,.
  \par
For $d\geq2$\,, it suffices to prove that
\begin{equation*}
  z_k-1 \in {<}\,\bigl\{ q_\bTheta : \bTheta \in \cU_d^{(k)}
  \bigr\}\,{>} \quad\text{for }k=1,\dots,d\,,
\end{equation*}
where the family $\cU_d^{(k)}$ is defined to be
\begin{equation*}
  \cU_d^{(k)} = \Bigl\{ \bTheta : d{\times}(d{-}1)\text{-submatrix of }
  \bigl(\bX_d^{(1)}
  \setlength{\unitlength}{1mm}\begin{picture}(2,3)
  \put(0,0){\dottedline{0.5}(1,-1)(1,3)}
  \end{picture} \bX_d^{(2)}\bigr)\,,\;\det \bigl(\bTheta
  \setlength{\unitlength}{1mm}\begin{picture}(2,3)
  \put(0,0){\dottedline{0.5}(1,-1)(1,3)}
  \end{picture} \be_k\bigr)= \pm1 \Bigr\} \,.
\end{equation*}
The proof is by induction on $d$\,. For $d=2$\,, we have
$\bigl(\bX_d^{(1)}
  \setlength{\unitlength}{1mm}\begin{picture}(2,3)
  \put(0,0){\dottedline{0.5}(1,-1)(1,3)}
  \end{picture} \bX_d^{(2)}\bigr) = \left(
  \begin{smallmatrix}1&0&1\\[0.5mm]0&1&1\end{smallmatrix} \right)$\,.
From this, we find for $k=1$
\begin{equation*}
  \cU_2^{(1)} = \Bigl\{
  \left(\begin{smallmatrix}0\\[1mm]1\end{smallmatrix}\right) \,{,}\,
  \left(\begin{smallmatrix}1\\[1mm]1\end{smallmatrix}\right) \Bigr\}
\end{equation*}
and thus the lemma claims that
\begin{equation*}
  z_1-1 \in {<}\,\mfrac{1+z_2}{2}\,{,}\,\mfrac{1+z_1\,z_2}{2} \,{>}
\end{equation*}
which follows from
\begin{equation}\label{eq:indstart}
  z_1-1 = 2\,z_1\,\mfrac{1+z_2}{2} - 2\,\mfrac{1+z_1\,z_2}{2} \,.
\end{equation}
The claim for $k=2$ follows by symmetry.
  \par
For the induction step, we consider first $k=d+1$ and write
\begin{equation}\label{eq:indstep}
  z_{d+1}-1 = -z_{d+1}\,(z_d-1)\,\mfrac{1+z_d}{2} +
  (z_d-1)\,\mfrac{1+z_d\,z_{d+1}}{2} + (z_{d+1}-1)\,\mfrac{1+z_d}{2}
  \,.
\end{equation}
By the induction hypothesis, we can write
\begin{equation}\label{eq:indhyp}
  z_d-1 = \sum_{\bTheta\in\cU_d^{(d)}}
  p_\bTheta(z_1,\dots,z_d)\,q_\bTheta(z_1,\dots,z_d)
\end{equation}
for certain polynomials $p_\bTheta$\,. For $\bTheta\in\cU_{d}^{(d)}$\,,
we define
\begin{equation*}
  \bTheta' = \begin{pmatrix}
                & \text{\footnotesize$0$} \\[-2.5mm]
    \,\bTheta\! & \text{\footnotesize$\vdots$} \\[-2mm]
                & \text{\footnotesize$0$} \\[-1mm]
                & \text{\footnotesize$1$} \\[-1mm]
    \text{\footnotesize$0\;\cdots\;0$}\!\!\!\!
                &\text{\footnotesize$0$}
  \end{pmatrix} \!,\quad
  \bTheta'' = \begin{pmatrix}
                & \text{\footnotesize$0$} \\[-2.5mm]
    \,\bTheta\! & \text{\footnotesize$\vdots$} \\[-2mm]
                & \text{\footnotesize$0$} \\[-1mm]
                & \text{\footnotesize$1$} \\[-1mm]
    \text{\footnotesize$0\;\cdots\;0$}\!\!\!\!
                &\text{\footnotesize$1$}
  \end{pmatrix} \!, \quad\text{and}\quad
  \bTheta''' = \begin{pmatrix}
                   & \text{\footnotesize$0$} \\[-2.5mm]
                   & \text{\footnotesize$\vdots$} \\[-2mm]
    \bP\,\bTheta\! & \text{\footnotesize$0$} \\[-1mm]
                   & \text{\footnotesize$1$} \\[-1mm]
                   & \text{\footnotesize$0$}
  \end{pmatrix} \!,
\end{equation*}
where
  \begin{equation*}
  \bP = \text{\footnotesize$\displaystyle
  \setlength{\arraycolsep}{2pt}
  \begin{pmatrix}
    1 &        & 0 & 0      \\[-2mm]
      & \ddots &   & \vdots \\[-1mm]
    0 &        & 1 & 0      \\[-1mm]
    0 & \cdots & 0 & 0      \\[-1mm]
    0 & \cdots & 0 & 1
  \end{pmatrix}$} \,.
\end{equation*}
It is easily seen that this yields $\bTheta'$\,, $\bTheta''$\,,
$\bTheta'''\in\cU_{d+1}^{(d+1)}$ with
\begin{align*}
  q_{\bTheta'}(z_1,\ldots,z_{d+1}) &=
    q_\bTheta(z_1,\ldots,z_d)\,\mfrac{1+z_d}{2} \,,\\
  q_{\bTheta''}(z_1,\ldots,z_{d+1}) &=
    q_\bTheta(z_1,\ldots,z_d)\,\mfrac{1+z_d\,z_{d+1}}{2} \,, \\
  \text{and}\qquad
    q_{\bTheta'''}(z_1,\ldots,z_{d+1}) &=
    q_\bTheta(z_1,\ldots,z_{d-1},z_{d+1})\,\mfrac{1+z_d}{2} \,.
\end{align*}
Using these identities, we obtain from~\eqref{eq:indhyp}
\begin{align*}
  -z_{d+1}\,(z_d-1)\,\mfrac{1+z_d}{2} &= \sum_{\bTheta\in\cU_d^{(d)}}
    -z_{d+1}\,p_\bTheta(z_1,\dots,z_d)\,q_{\bTheta'}(z_1,\dots,z_{d+1})
    \,,\\
  (z_d-1)\,\mfrac{1+z_d\,z_{d+1}}{2} &= \sum_{\bTheta\in\cU_d^{(d)}}
    p_\bTheta(z_1,\dots,z_d)\,q_{\bTheta''}(z_1,\dots,z_{d+1}) \,,
\intertext{and, by replacing $z_d$ by $z_{d+1}$ in~\eqref{eq:indhyp},}
  (z_{d+1}-1)\,\mfrac{1+z_d}{2} &= \sum_{\bTheta\in\cU_d^{(d)}}
    p_\bTheta(z_1,\dots,z_{d-1},z_{d+1})\,
    q_{\bTheta'''}(z_1,\dots,z_{d+1}) \,.
\end{align*}
Substituting these three identities into~\eqref{eq:indstep} yields a
representation for $z_{d+1}-1$ of the desired form. A corresponding
representation for $z_k-1$ with $k\leq d$ can be obtained by a cyclic
permutation of the indices.

This completes the induction.
\end{proof}

We would like to emphasize again that a very important consequence of
the following theorem, for subdivision schemes, is stated as the main
result of this section in Theorem~A.

\begin{theorem}\label{thm:JBoxgen}
The ideal $\cI$ is generated by the elements $q_\bTheta$\,,
where $\bTheta$ are the unimodular $d{\times}d$-submatrices of
$\bigl(\bX_d^{(1)}
\setlength{\unitlength}{1mm}\begin{picture}(2,3)
\put(0,0){\dottedline{0.5}(1,-1)(1,3)}
\end{picture} \bX_d^{(2)}\bigr)$ from~\eqref{eq:directions}.
\end{theorem}

\begin{proof}
One of the inclusions follows from Proposition~\ref{prp:nonunimodcase}.
In the light of~\eqref{eq:genideal}, the other inclusion is an
immediate consequence of Lemma~\ref{lem:Jgensingle} and the fact that
\begin{equation*}
  (z_1+1)(z_2+1)\ldots(z_d+1) = 2^d\,q_{\bTheta}(\bz)
  \quad\text{for}\quad \bTheta = \bI_d = \bX_d^{(1)} \,.
\end{equation*}
\end{proof}

It is worth noting that the condition
$\be_k \in\bTheta$ in
Lemma~\ref{lem:Jgensingle} is equivalent to
the divisibility of $q_\bTheta$ by
$1{+}z_k $\,, which seems natural since
$1{+}z_k$ obviously divides $1{-}z_k^2$\,. Furthermore, the fact that each matrix
$\bTheta$ in Theorem~\ref{thm:JBoxgen} contains a
 standard
unit vector as its column is a consequence of the following result.

\begin{lemma}\label{lem:evenDetinX2}
Any $d{\times}d$-submatrix $\bTheta$ of $\bX_d^{(2)}$ satisfies
\begin{equation*}
  \det \bTheta \equiv 0 \;(\operatorname{mod}\; 2)\;.
\end{equation*}
\end{lemma}

\begin{proof}
A $d{\times}d$-submatrix $\bTheta$ of $\bX_d^{(2)}$ contains in
each column exactly two entries equal to one (and the others are
equal to  zero). For the linear operator $L$ introduced in the
proof of Proposition~\ref{prp:nonunimodcase}, this implies
\begin{equation*}
  L(\bone^T) = \bone^T\bTheta = 2\cdot\bone^T \equiv \bnull^T
  \;(\operatorname{mod}\; 2)\,,
\end{equation*}
\textit{i.\,e.}, $\bone^T \in \ker (\overline{L})$\,. From this,
we may conclude that $\det \overline{L} = \overline{0}$ in $\ZZ_2$
and, thus, $\det\bTheta = \det L \equiv 0 \;(\operatorname{mod}\;
2)$\,.
\end{proof}

\begin{examples}
We list the generators $q_\bTheta$ of $\cI$ for low-dimensional
cases. {\setlength{\leftmargini}{5mm}
   \begin{itemize}
\item
For $d=1$\,, we have $\bX_1=(1)$\,, so the only submatrix is
$\bTheta=\bX_1$ with
\begin{equation*}
  q_\bTheta(z) = \mfrac{1+z}{2}\,,
\end{equation*}
and $\cI$ is the principal ideal generated by this function.
\item
For $d=2$, we have
\begin{equation*}
  \bX_2 =
  \setlength{\unitlength}{1mm}\begin{picture}(0,0)
  \put(0,0){\dottedline{0.5}(12,-3)(12,5)}
  \end{picture}
  \begin{pmatrix} 1 & 0 & 1 \\ 0 & 1 & 1 \end{pmatrix}
\end{equation*}
which yields the following three generators for~$\cI$\,:
\begin{alignat*}{3}
  \bTheta_1 &= \begin{pmatrix}1&0\\0&1\end{pmatrix}
    \qquad\text{with}\qquad
    & q_{\bTheta_1}(z) &= \mfrac{1+z_1}{2}\;\mfrac{1+z_2}{2}\;,\\
  \bTheta_2 &= \begin{pmatrix}1&1\\0&1\end{pmatrix}
    \qquad\text{with}\qquad
    & q_{\bTheta_2}(z) &= \mfrac{1+z_1}{2}\;\mfrac{1+z_1\,z_2}{2}\;,\\
  \bTheta_3 &= \begin{pmatrix}0&1\\1&1\end{pmatrix}
    \qquad\text{with}\qquad
    & q_{\bTheta_3}(z) &= \mfrac{1+z_2}{2}\;\mfrac{1+z_1\,z_2}{2}\;,
\end{alignat*}
as stated in~\eqref{eq:Jd2Box}. These three functions generate $\cI$
minimally in the sense that no two of them generate all of $\cI$\,.
\item
In the case $d=3$, we have
\begin{equation*}
  \bX_3 =
  \setlength{\unitlength}{1mm}\begin{picture}(0,0)
  \put(0,0){\dottedline{0.5}(18,-5)(18,7)}
  \put(0,0){\dottedline{0.5}(34,-5)(34,7)}
  \end{picture}
  \begin{pmatrix}
    1 & 0 & 0 & 1 & 1 & 0 & 1 \\
    0 & 1 & 0 & 1 & 0 & 1 & 1 \\
    0 & 0 & 1 & 0 & 1 & 1 & 1
  \end{pmatrix}.
\end{equation*}
The submatrices $\bTheta$ with even determinants are obtained by
selecting the columns $(1,2,4)$, $(1,3,5)$, $(1,6,7)$, $(2,3,6)$,
$(2,5,7)$, and $(3,4,7)$,  $\det \bTheta = 0$, and $(4,5,6)$,
  $\det \bTheta = -2$ (compare Lemma~\ref{lem:evenDetinX2}). The
remaining $\binom{7}{3} - 7 = 28$ submatrices have determinant $\pm1$
and thus describe elements of~$\cI$\,. Restricting ourselves to
$\bigl(\bX_d^{(1)} \setlength{\unitlength}{1mm}\begin{picture}(2,3)
\put(0,0){\dottedline{0.5}(1,-1)(1,3)}
\end{picture} \bX_d^{(2)}\bigr)$ as stated in Theorem~\ref{thm:JBoxgen}
yields $\binom{6}{3} - 4 = 16$ elements generating the ideal $\cI$\,.
However, this system is highly redundant: using an algebraic
manipulation program, we found that it contains $64$ subsets of $10$
elements each which generate $\cI$ minimally (in the sense that in each
case, omitting any one of the $10$ elements does not yield a set of
generators anymore).
\end{itemize}}
\end{examples}

\subsection*{Geometric interpretation and modification.}

The columns of $\bX_d$ describe the edges and the various
diagonals of the $d$-hypercube. By
Proposition~\ref{prp:nonunimodcase}, the polynomial $q_\bTheta$ is
an element of $\cI$ if and only if the columns of $\bTheta$ span a
parallelepiped with an odd $d$-volume. Theorem~\ref{thm:JBoxgen}
states that in order to generate $\cI$\,, it suffices to restrict
ourselves to the edges $\be_k$ and the $2$-surface-diagonals $\be_j+\be_k$, $j\ne k$, of the
hypercube only, and to consider only parallelepipeds with
$d$-volume equal to one. Lemma~\ref{lem:evenDetinX2} shows that
these parallelepipeds all have at least one unit vector as an
edge.

\smallskip
Proposition~\ref{prp:Jd2Boxmod} leads to an alternative approach. For $d=2$,
the three vectors
$\be_1 = \left(\begin{smallmatrix}1\\[0.5mm]0\end{smallmatrix}\right)$,
$\be_2 = \left(\begin{smallmatrix}0\\[0.5mm]1\end{smallmatrix}\right)$, and
$\be_1-\be_2 = \left(\begin{smallmatrix}1\\[0.5mm]-1\end{smallmatrix}\right)$
can also be seen as the edges of the $2$-simplex. Since for $d>2$ we only
need to use the diagonals of the
$2$-dimensional surfaces of the $d$-hypercube, it seems natural to consider
the edges of the $d$-simplex instead.

To this end, we describe a different family of polynomials. For
any $d{\times}d$-submatrix $\bTheta$ of $\bigl(\bX_d^{(1)}
\setlength{\unitlength}{1mm}\begin{picture}(2,3)
\put(0,0){\dottedline{0.5}(1,-1)(1,3)}
\end{picture} \bX_d^{(2)}\bigr)$\,, define
\begin{equation*}
  \widetilde{q}_\bTheta(\bz) = \prod_{\btheta\in\bTheta}
  \widetilde{r}_\btheta(\bz) \quad\text{with}\quad
  \widetilde{r}_\btheta(\bz) = \begin{cases}
    \mfrac{1+z_k}{2}\,, & \text{if } \btheta=\be_k\,, \\[2mm]
    \mfrac{z_j+z_k}{2}\,, & \text{if } \btheta=\be_j+\be_k\,.
  \end{cases}
\end{equation*}
Geometrically, this amounts for $d=2$ to replacing the vector
$\left(\begin{smallmatrix}1\\[0.5mm]1\end{smallmatrix}\right)$ by
$\left(\begin{smallmatrix}1\\[0.5mm]-1\end{smallmatrix}\right)$ as
third direction in the three-directional grid. Equivalently, this can
be seen as a reflection of the grid about one of the coordinate axes
(compare Remark~\ref{rem:gridrefl}).

Algebraically, this has the advantage that each factor
$\widetilde{r}_\btheta$ is linear and therefore that all the
$\widetilde{q}_\bTheta$ share the same total degree~$d$\,. These
functions still satisfy the analogues of
Proposition~\ref{prp:nonunimodcase}, Lemma~\ref{lem:Jgensingle}, and,
consequently, of Theorem~\ref{thm:JBoxgen} and
Theorem~A.

\begin{proposition}\label{prp:nonunimodcasemod}
For any $d{\times}d$-submatrix $\bTheta$ of $\bigl(\bX_d^{(1)}
\setlength{\unitlength}{1mm}\begin{picture}(2,3)
\put(0,0){\dottedline{0.5}(1,-1)(1,3)}
\end{picture} \bX_d^{(2)}\bigr)$ in \eqref{eq:directions}, we have
\begin{equation*}
  \widetilde{q}_\bTheta \in \cI \quad \Longleftrightarrow \quad
  \det \bTheta \equiv 1 \;(\operatorname{mod}\; 2)\;.
\end{equation*}
\end{proposition}

\begin{proof}
For $\beps=(\varepsilon_1, \dots, \varepsilon_d)\in\rZ'=\{{-}1\,{,}\,{+}1\}^d\setminus\{\bone\}$\,, we have
\begin{equation*}
  \mfrac{\varepsilon_j+\varepsilon_k}{2} =
  \varepsilon_k\,\mfrac{1+\varepsilon_j/\varepsilon_k}{2} =
  \pm\mfrac{1+\varepsilon_j\,\varepsilon_k}{2}, \quad j,k=1,\dots,d,
\end{equation*}
and, therefore, $\widetilde{q}_\bTheta(\beps) = \pm q_\bTheta(\beps)$\,.
Consequently,
\begin{equation*}
  \widetilde{q}_\bTheta \in \cI \iff
  \widetilde{q}_\bTheta\big|_{\rZ'}=0 \iff
  q_\bTheta\big|_{\rZ'}=0 \iff
  q_\bTheta \in \cI\,.
\end{equation*}
The claim follows by Proposition~\ref{prp:nonunimodcase}.
\end{proof}

\begin{lemma}\label{lem:Jgensinglemod}
For $k=1,\dots,d$\,, we have
\begin{equation*}
  1-z_k^2 \;\in\; {<}\,\{ \widetilde{q}_\bTheta : \bTheta\subseteq
  \bigl(\bX_d^{(1)}
  \setlength{\unitlength}{1mm}\begin{picture}(2,3)
  \put(0,0){\dottedline{0.5}(1,-1)(1,3)}
  \end{picture} \bX_d^{(2)}\bigr)\,,\;\det \bTheta
  = \pm1\,,\; \be_k \in\bTheta \}\,{>} \,.
\end{equation*}
\end{lemma}

\begin{proof}
The proof is \textit{mutatis mutandis} the same as that of
Lemma~\ref{lem:Jgensingle}, except that~\eqref{eq:indstart} becomes
\begin{equation*}
  z_1-1 = -2\,\mfrac{1+z_2}{2} + 2\,\mfrac{z_1+z_2}{2}
\end{equation*}
and~\eqref{eq:indstep} becomes
\begin{equation*}
  z_{d+1}-1 = (z_d-1)\,\mfrac{1+z_d}{2} -
  (z_d-1)\,\mfrac{z_d+z_{d+1}}{2} + (z_{d+1}-1)\,\mfrac{1+z_d}{2}
  \,.
\end{equation*}
\end{proof}

\begin{theorem}\label{thm:JBoxgenmod}
The ideal $\cI$ is generated by the elements
$\widetilde{q}_\bTheta$\,, where $\bTheta$ are the unimodular
$d{\times}d$-submatrices of $\bigl(\bX_d^{(1)}
\setlength{\unitlength}{1mm}\begin{picture}(2,3)
\put(0,0){\dottedline{0.5}(1,-1)(1,3)}
\end{picture} \bX_d^{(2)}\bigr)$ from~\eqref{eq:directions}.
\end{theorem}

\begin{proof}
As in the proof of Theorem~\ref{thm:JBoxgen}, using
Lemma~\ref{lem:Jgensinglemod} instead of Lemma~\ref{lem:Jgensingle}.
\end{proof}

\medskip
The interpretation of Theorem~\ref{thm:JBoxgenmod}   as a property of
subdivision schemes
leads to the following modification of Theorem~A:

\medskip
\noindent {\bf Theorem $\widetilde{\bf A}$.} {\it The mask symbol
of any convergent $d$-variate subdivision scheme $\cS_\ra$ can be
written in the form
\begin{equation*}
  \ra(\bz) = \sum_{\bTheta} \lambda_\bTheta\; \sigma_\bTheta(\bz)\;
  2^d\,\widetilde{q}_\bTheta(\bz)\;,
\end{equation*}
where $\sigma_\bTheta(\bz)$ are Laurent polynomials satisfying
$\sigma_\bTheta(\bone)=1$\,, and $\lambda_\bTheta$ are real numbers
subject to $\sum_{\bTheta} \lambda_\bTheta = 1$. The sum runs over
all unimodular $d{\times}d$-submatrices $\bTheta$ of $\bigl(\bX_d^{(1)}
\setlength{\unitlength}{1mm}\begin{picture}(2,3)
\put(0,0){\dottedline{0.5}(1,-1)(1,3)}
\end{picture} \bX_d^{(2)}\bigr)$ from
\eqref{eq:directions}.}

\section{Zero condition of higher orders and powers of $\cI$ for $d=2$} 

In this section, we describe families of generators for the
ideals~$\cI^k$\,, for $k>1$, and the implications on the representation of
mask symbols from Theorem~A. The case $d=2$ is of special interest in
subdivision. Therefore,
most of this section is devoted to the study of this case, see in particular
Theorem~C below
and the remarks at the end of Section~3.1.



We start with a simple observation from ideal theory, namely, the fact that
the product  $\cI_1\cdot\cI_2$ of two ideals
$\cI_1 = {<}\,\{a_j:j=1,\dots,n\}\,{>}$ and
$\cI_2 = {<}\,\{b_k:k=1,\dots,m\}\,{>}$ in a   ring is generated by the
pointwise products of the corresponding generating sets
\begin{equation} \label{eq:Iprodgen}
 \cI_1\cdot\cI_2 = {<}\, \{a_j\,b_k : j=1,\dots,n,
  k=1,\dots,m\}\,{>}\;.
\end{equation}
Applying this to the construction of generators for $\cI^k$, we get the
following immediate generalization of
Theorem~A.

\medskip
\noindent
{\bf Theorem B.} {\it
A convergent $d$-variate subdivision scheme $\cS_\ra$ satisfies
the condition~$\rZ_k$ if and only if its mask symbol can be
written in the form
\begin{equation*}
  \ra(\bz) = \sum_{j} \lambda_j \; \sigma_j(\bz) \; 2^d\,Q_j(\bz)\,,
\end{equation*}
where $Q_j(\bz)$ are $k$-fold products of Laurent polynomials
$q_\bTheta$ with unimodular $d{\times}d$-submatrices $\bTheta$ of $X_d$
from \eqref{eq:directions}, $\sigma_j(\bz)$ are Laurent
polynomials normalized by $\sigma_j(\bone)=1$\,, and $\lambda_j$ are
real numbers subject to $\sum_j \lambda_j = 1$\,.}

\medskip
Note that the normalizations assumed in Theorem~B imply
$\ra(\bone) = 2^d$\,.

\bigskip
In the bivariate case, we can be much more specific and show that the
generators for $\cI^k$ are the mask symbols of certain
three-directional box splines. The latter have been studied thoroughly
in \cite{BH2,BHR,Chui,LaiSch}, and their mask symbols have the form
$B_{\alpha,\beta,\gamma} = 4\,B_{\alpha,\beta,\gamma}^{\#}$\,, where
\begin{equation} \label{eq:threedirsymbols}
  B_{\alpha,\beta,\gamma}^{\#}(z_1,z_2) =
  \Bigl(\mfrac{1+z_1}{2}\Bigr)^\alpha \,
  \Bigl(\mfrac{1+z_2}{2}\Bigr)^\beta \,
  \Bigl(\mfrac{1+z_1\,z_2}{2}\Bigr)^\gamma
, \quad \alpha,\beta,\gamma\in\NN_0,
\end{equation}
are the normalized box spline symbols satisfying
$B_{\alpha,\beta,\gamma}^{\#}(1,1) = 1$\,.
The three indices
$\alpha,\beta,\gamma$ correspond
to the multiplicities of the three vectors $\be_1$, $\be_2$ and $\be_1+\be_2$.
In this
notation, \eqref{eq:Jd2Box} becomes
\begin{equation}\label{eq:J1Boxgen}
  \cI = {<}\, B_{1,1,0}^{\#} \,, B_{1,0,1}^{\#}   \,, B_{0,1,1}^{\#} \,{>}
  \,.
\end{equation}
Note that the family of these (normalized) box spline symbols is partially
ordered
and closed under multiplication, since
\begin{equation*}
  B_{\alpha,\beta,\gamma}^{\#} \cdot B_{\alpha',\beta',\gamma'}^{\#} =
  B_{\alpha+\alpha',\beta+\beta',\gamma+\gamma'}^{\#}\,.
\end{equation*}
Furthermore, they satisfy the following relation that we  need later on.

\begin{lemma}\label{lem:genred}
For any given triple $(\alpha,\beta,\gamma) \in\NN_0^3$\,, the ideal generated by
the three symbols $B_{\alpha+1,\beta,\gamma}^{\#}$\,,
$B_{\alpha,\beta+1,\gamma}^{\#}$\,, and
$B_{\alpha,\beta,\gamma+1}^{\#}$
is the principal ideal generated by $B_{\alpha,\beta,\gamma}^{\#}$\,.
\end{lemma}

\begin{proof}
Since each of the symbols $B_{\alpha+1,\beta,\gamma}^{\#}$\,,
$B_{\alpha,\beta+1,\gamma}^{\#}$\,, and
$B_{\alpha,\beta,\gamma+1}^{\#}$ is a multiple of
$B_{\alpha,\beta,\gamma}^{\#}$\,, we only have to show that the latter
can be generated from the former three. To this end, we make use of the
identity
\begin{equation}\label{eq:genred}
  \mfrac{1}{2}\,(1-z_2)\,B_{1,0,0}^{\#}(z_1,z_2) +
  \mfrac{1}{2}\,(1-z_1)\,B_{0,1,0}^{\#}(z_1,z_2) +
  B_{0,0,1}^{\#}(z_1,z_2) = 1 \,.
\end{equation}
Multiplying both sides by $B_{\alpha,\beta,\gamma}^{\#}$ proves the
lemma.
\end{proof}

In the light of~\eqref{eq:J1Boxgen}, it is natural to  expect that the higher
powers of the ideal~$\cI$ are generated by box spline symbols of higher
order, as the following result shows.

\begin{theorem}\label{thm:twoDJkgen}
In the bivariate case, the $k$-th power $\cI^k$, $k \in \NN$, of
the ideal $\cI$ is generated by the set of three-directional box
spline symbols
\begin{equation*}
  \rI_k := \Bigl\{\; B^\#_{\beta,\beta,\alpha},
  B^\#_{\beta,\alpha,\beta},
  B^\#_{\alpha,\beta,\beta}\;:\;
  \alpha=0,1,\dots,\Bigl\lfloor \mfrac{k}{2}
  \Bigr\rfloor,\, \ \beta=k\,{-}\,\alpha\;\Bigr\}\;.
\end{equation*}
\end{theorem}

\begin{proof}
The proof is  by induction on $k$\,. For $k=1$, the claim is
just the identity~\eqref{eq:J1Boxgen}.
  \par
For the induction step, we write $\cI^{k+1} = \cI\cdot\cI^k$ and
apply~\eqref{eq:Iprodgen}. Using once more~\eqref{eq:J1Boxgen} and the
induction hypothesis $\cI^k = {<}\,\rI_k\,{>}$ yields $\cI^{k+1} =
{<}\,\rI'_{k+1}\,{>}$ with
\begin{equation*}
  \rI'_{k+1} = \left\{ \begin{array}{l}
    B_{\beta+1,\beta+1,\alpha}^{\#}\,,\,
    B_{\beta+1,\alpha+1,\beta}^{\#}\,,\,
    B_{\alpha+1,\beta+1,\beta}^{\#}, \\
    B_{\beta+1,\beta,\alpha+1}^{\#}\,,\,
    B_{\beta+1,\alpha,\beta+1}^{\#}\,,\,
    B_{\alpha+1,\beta,\beta+1}^{\#}, \\
    B_{\beta,\beta+1,\alpha+1}^{\#}\,,\,
    B_{\beta,\alpha+1,\beta+1}^{\#}\,,\,
    B_{\alpha,\beta+1,\beta+1}^{\#},
  \end{array} \alpha = 0,1,\dots,\Bigl\lfloor \mfrac{k}{2}
  \Bigr\rfloor,\,\beta=k\,{-}\,\alpha
  \right\} \,.
\end{equation*}
So we need to show that ${<}\,\rI'_{k+1}\,{>} =
{<}\,\rI_{k+1}\,{>}$ where
\begin{equation*}
  \rI_{k+1} = \Bigl\{\; B^\#_{\delta,\delta,\gamma},
  B^\#_{\delta,\gamma,\delta},
  B^\#_{\gamma,\delta,\delta}\;:\;
  \gamma=0,1,\dots,\Bigl\lfloor \mfrac{k{+}1}{2}
  \Bigr\rfloor,\,\delta=k\,{+}\,1\,{-}\,\gamma\;\Bigr\}\;.
\end{equation*}
To this end, it suffices to show that both
$\rI_{k+1}\subseteq{<}\,\rI'_{k+1}\,{>}$ and
$\rI'_{k+1} \subseteq {<}\,\rI_{k+1}\,{>}$\,.
  \par
Firstly, note that the elements of $\rI_{k+1}$ are the diagonal
elements in the list $\rI'_{k+1}$ with $\gamma=\alpha$ and
thus $\delta=\beta+1$\,, with the only exception $\gamma=\Bigl\lfloor
\mfrac{k{+}1}{2}\Bigr\rfloor > \Bigl\lfloor\mfrac{k}{2}\Bigr\rfloor$\,.
This can only happen if $k=2\,\ell\,{+}\,1$ is odd, and then
$\gamma=\ell\,{+}\,1$\,. For this value of $\gamma$\,, the list
$\rI_{k+1}$ contains only one element, \textit{viz.},
$B^\#_{\ell+1,\ell+1,\ell+1}$\,. By Lemma~\ref{lem:genred}, this
is generated by the three elements above the main diagonal in
$\rI'_{k+1}$ with $k=2\,\ell\,{+}\,1$ and
$\alpha=\ell\,,\;\beta=\ell{+}1$\,.
  \par
Conversely, the diagonal elements of $\rI'_{k+1}$ are all
listed in $\rI_{k+1}$\,. The indices of any non-diagonal element form a
permutation of the triple $(\beta,\beta{+}1,\alpha{+}1)$\,. The
associated $B_{\beta,\beta+1,\alpha{+}1}^{\#}$ is a multiple of
$B_{\beta,\beta,\alpha{+}1}^{\#}$ which appears together with all index
permutations in $\rI_{k+1}$ for $\gamma=\alpha{+}1\,,\;\delta=\beta$\,,
except for the case $\alpha = \Bigl\lfloor\mfrac{k}{2}\Bigr\rfloor =
\Bigl\lfloor\mfrac{k{+}1}{2}\Bigr\rfloor$\,. This can only happen if
$k=2\,\ell$ is even, and then $\alpha=\ell$ and
$(\beta,\beta{+}1,\alpha{+}1)=(\ell,\ell{+}1,\ell{+}1)$\,. But the
associated elements appear in $\rI_{k+1}$ for
$\gamma=\ell\,,\;\delta=\ell{+}1$\,.
\end{proof}

As in Section~\ref{sec:zeroconideal}, Theorem~\ref{thm:twoDJkgen} has
an immediate consequence, Theorem~C, for bivariate mask symbols satisfying
the higher order zero conditions. This result is of great importance for
studying the properties of existing subdivision schemes and also as a starting
point for the construction of new schemes.

\medskip
\noindent
{\bf Theorem C.} {\it
A convergent bivariate subdivision scheme $\cS_\ra$ satisfies the
condition~$\rZ_k$ if and only if its mask symbol can be
written in the form
\begin{equation}\label{eq:bivgen}
  \ra(\bz) = \sum_{B^\#_{\alpha,\beta,\gamma}\in\rI_k}
  \lambda_{\alpha,\beta,\gamma} \;
  \sigma_{\alpha,\beta,\gamma}(\bz) \;
  4\,B^\#_{\alpha,\beta,\gamma}(\bz)\,,
\end{equation}
where $\sum \lambda_{\alpha,\beta,\gamma} = 1$\,, and the
$\sigma_{\alpha,\beta,\gamma}(\bz)$ are Laurent polynomials normalized
by the condition $\sigma_{\alpha,\beta,\gamma}(\bone)=1$\,.}

\begin{examples}\label{exas:twoDJkgen}
We illustrate the result  of Theorem~\ref{thm:twoDJkgen}
 by explicitly listing the generators for
small values of $k$\,.
\begin{align*}
  \rI_1 &= \{ \, B^\#_{1,1,0} \,,\, B^\#_{1,0,1} \,,\, B^\#_{0,1,1}
    \, \} \,, \\[1mm]
  \rI_2 &=\{ \, B^\#_{2,2,0} \,,\, B^\#_{2,0,2} \,,\, B^\#_{0,2,2}
    \,,\, B^\#_{1,1,1} \, \} \,, \\[1mm]
  \rI_3 &= \{ \, B^\#_{3,3,0} \,,\, B^\#_{3,0,3} \,,\, B^\#_{0,3,3}
    \,,\, B^\#_{2,2,1} \,,\, B^\#_{2,1,2} \,,\, B^\#_{1,2,2} \, \} \,,
    \\[1mm]
  \rI_4 &= \{ \, B^\#_{4,4,0} \,,\, B^\#_{4,0,4} \,,\, B^\#_{0,4,4}
    \,,\, B^\#_{3,3,1} \,,\, B^\#_{3,1,3} \,,\, B^\#_{1,3,3} \,,\,
    B^\#_{2,2,2} \, \} \,. \\[1mm]
\end{align*}
\end{examples}
\subsection{Further properties. \label{ssec:prop}}

To be able to show some further properties of the box spline symbols,
 implying the
corresponding properties of the associated subdivision schemes,
we need the following auxiliary result. As usual, we write
\begin{equation*}
  D^{(n,m)} =
  \mfrac{\partial^{n+m}}{\partial z_1^{\;n}\;\partial z_2^{\;m}}
, \quad n,m \in \NN_0,
\end{equation*}
for mixed partial differential operators.

\begin{lemma}\label{lem:BoxspliDer}
The partial derivatives of the box spline symbol
$B^{\#}_{\alpha,\beta,\gamma}$ are given by
\begin{align}
  \Bigl(D^{(n,m)}
    \, B^{\#}_{\alpha,\beta,\gamma}\Bigr) &(z_1,z_2) = \notag\\
  = \sum_{\ell=0}^\gamma \mbinom{\gamma}{\ell}
    & \biggl( \sum_{i=0}^n \mfrac{n!}{2^n}
    \mbinom{\alpha{+}\ell}{n{-}i}
    \Bigl( \mfrac{1{+}z_1}{2} \Bigr)^{\alpha+\ell-(n-i)}
    \mbinom{\gamma{-}\ell}{i}
    \Bigl( \mfrac{z_1{-}1}{2} \Bigr)^{\gamma-\ell-i} \biggr)
    \label{eq:BoxspliDer}\\
  &\times \biggl( \sum_{j=0}^m \mfrac{m!}{2^m}
    \mbinom{\beta{+}\ell}{m{-}j}
    \Bigl( \mfrac{1{+}z_2}{2} \Bigr)^{\beta+\ell-(m-j)}
    \mbinom{\gamma{-}\ell}{j}
    \Bigl( \mfrac{z_2{-}1}{2} \Bigr)^{\gamma-\ell-j} \biggr)\,. \notag
\end{align}
\end{lemma}

\begin{proof}
The identity
\begin{equation*}
   \mfrac{1+z_1\,z_2}{2} = \mfrac{1+z_1}{2}\,\mfrac{1+z_2}{2} +
   \mfrac{z_1-1}{2}\,\mfrac{z_2-1}{2}
\end{equation*}
yields
\begin{align*}
  B^{\#}_{\alpha,\beta,\gamma} (z_1,z_2) &=
    \Bigl( \mfrac{1{+}z_1}{2} \Bigr)^\alpha
    \Bigl( \mfrac{1{+}z_2}{2} \Bigr)^\beta
    \Bigl( \mfrac{1+z_1}{2}\,\mfrac{1+z_2}{2} +
    \mfrac{z_1-1}{2}\,\mfrac{z_2-1}{2} \Bigr)^\gamma \\
  &= \sum_{\ell=0}^\gamma \mbinom{\gamma}{\ell}
    \Bigl( \mfrac{1{+}z_1}{2} \Bigr)^{\alpha+\ell}
    \Bigl( \mfrac{1{+}z_2}{2} \Bigr)^{\beta+\ell}
    \Bigl( \mfrac{z_1{-}1}{2} \Bigr)^{\gamma-\ell}
    \Bigl( \mfrac{z_2{-}1}{2} \Bigr)^{\gamma-\ell} \,.
\end{align*}
From this we obtain with the Leibniz formula
\begin{align*}
  \Bigl(D^{(n,m)}
    \, B^{\#}_{\alpha,\beta,\gamma}\Bigr) & (z_1,z_2) = \\
  = \sum_{\ell=0}^\gamma \mbinom{\gamma}{\ell}
    & \mfrac{\partial^n}{\partial z_1^{\;n}}
    \biggl( \Bigl( \mfrac{1{+}z_1}{2} \Bigr)^{\alpha+\ell}
    \Bigl( \mfrac{z_1{-}1}{2} \Bigr)^{\gamma-\ell} \biggr)
    \mfrac{\partial^m}{\partial z_2^{\;m}}
    \biggl( \Bigl( \mfrac{1{+}z_2}{2} \Bigr)^{\beta+\ell}
    \Bigl( \mfrac{z_2{-}1}{2} \Bigr)^{\gamma-\ell} \biggr) \\
  = \sum_{\ell=0}^\gamma \mbinom{\gamma}{\ell}
    & \biggl( \sum_{i=0}^n \mbinom{n}{i}
    \mfrac{(\alpha{+}\ell)!}{(\alpha{+}\ell{-}(n{-}i))!}
    \mfrac{(1{+}z_1)^{\alpha+\ell-(n-i)}}{2^{\alpha+\ell}_{\nix}}
    \mfrac{(\gamma{-}\ell)!}{(\gamma{-}\ell{-}i)!}
    \mfrac{(z_1{-}1)^{\gamma-\ell-i}}{2^{\gamma-\ell}_{\nix}} \biggr) \\
  &\times \biggl( \sum_{j=0}^m \mbinom{m}{j}
    \mfrac{(\beta{+}\ell)!}{(\beta{+}\ell{-}(m{-}j))!}
    \mfrac{(1{+}z_2)^{\beta+\ell-(m-j)}}{2^{\beta+\ell}_{\nix}}
    \mfrac{(\gamma{-}\ell)!}{(\gamma{-}\ell{-}j)!}
    \mfrac{(z_2{-}1)^{\gamma-\ell-j}}{2^{\gamma-\ell}_{\nix}} \biggr) \\
  = \sum_{\ell=0}^\gamma \mbinom{\gamma}{\ell}
    & \biggl( \sum_{i=0}^n \mfrac{n!}{2^n}
    \mbinom{\alpha{+}\ell}{n{-}i}
    \Bigl( \mfrac{1{+}z_1}{2} \Bigr)^{\alpha+\ell-(n-i)}
    \mbinom{\gamma{-}\ell}{i}
    \Bigl( \mfrac{z_1{-}1}{2} \Bigr)^{\gamma-\ell-i} \biggr) \\
  &\times \biggl( \sum_{j=0}^m \mfrac{m!}{2^m}
    \mbinom{\beta{+}\ell}{m{-}j}
    \Bigl( \mfrac{1{+}z_2}{2} \Bigr)^{\beta+\ell-(m-j)}
    \mbinom{\gamma{-}\ell}{j}
    \Bigl( \mfrac{z_2{-}1}{2} \Bigr)^{\gamma-\ell-j} \biggr)
\end{align*}
as claimed.
\end{proof}

Together with Theorem~\ref{thm:twoDJkgen}, this allows us to determine
the maximal order of sum rules satisfied by a three-directional
box spline symbol.

\begin{proposition}\label{prp:maxJk}
For any triple $(\alpha,\beta,\gamma)\in\NN_0^3$\,, the maximal $k$
such that $B^{\#}_{\alpha,\beta,\gamma}\in\cI^k$ is given by
\begin{equation*}
  k = \alpha+\beta+\gamma-\max\{\alpha,\beta,\gamma\}\,.
\end{equation*}
\end{proposition}

\begin{proof}
Let $(\alpha,\beta,\gamma)=\pi(\alpha',\beta',\gamma')$ where $\pi$ is
a permutation such that $\alpha'\leq\beta'\leq\gamma'$\,. It follows
from Theorem~\ref{thm:twoDJkgen} that
\begin{equation*}
  B^{\#}_{\alpha,\beta,\gamma} =
  B^{\#}_{\pi(\alpha',\beta',\gamma')} =
  B^{\#}_{\pi(0,0,\gamma'-\beta')} \cdot
  B^{\#}_{\pi(\alpha',\beta',\beta')} \in \cI^k
\end{equation*}
for $k = \alpha'+\beta' = \alpha+\beta+\gamma -
\max\{\alpha,\beta,\gamma\}$\,.
  \par
On the other hand, applying~\eqref{eq:BoxspliDer} yields
\begin{alignat*}{3}
  &\Bigl(D^{(\alpha,\beta)}
    \, B^{\#}_{\alpha,\beta,\gamma}\Bigr) ({-}1\,{,}\,{-}1) &&=&
    \mfrac{\alpha!\,\beta!}{2^{\alpha+\beta}_{\nix}} &\not= 0 \,,\\
  &\Bigl(D^{(\alpha,\gamma)}
    \, B^{\#}_{\alpha,\beta,\gamma}\Bigr) ({-}1\,{,}\,1) &&=& \;
    (-1)^\gamma\mfrac{\alpha!\,\gamma!}{2^{\alpha+\gamma}_{\nix}}
    &\not= 0 \,,\\
  \text{and}\quad
  &\Bigl(D^{(\gamma,\beta)}
    \, B^{\#}_{\alpha,\beta,\gamma}\Bigr) (1\,{,}\,{-}1) &&=&
    (-1)^\gamma\mfrac{\beta!\,\gamma!}{2^{\beta+\gamma}_{\nix}}
    &\not= 0 \,.
\end{alignat*}
This shows that $B^{\#}_{\alpha,\beta,\gamma}\notin\cI^{k+1}$ for
\begin{equation*}
  k=\min\{\alpha{+}\beta\,,\,\alpha{+}\gamma\,,\,\beta{+}\gamma\} =
  \alpha+\beta+\gamma-\max\{\alpha,\beta,\gamma\}\,,
\end{equation*}
and this completes the proof.
\end{proof}

The sets of generators described in Theorem~\ref{thm:twoDJkgen} are
minimal generating sets, as the following shows.

\begin{proposition}\label{prp:mingen}
The set $\rI_k$ generating $\cI^k$ is minimal in the sense that for any
$B^\#_{\alpha,\beta,\gamma}\in\rI_k$\,, the reduced set
$\rI_k\setminus\{B^\#_{\alpha,\beta,\gamma}\}$ does no longer generate
the ideal~$\cI^k$\,.
\end{proposition}

\begin{proof}
We begin with an element of the form $B^\#_{\alpha,\beta,\beta}
\in\rI_k$\,. We find that
\begin{equation}\label{eq:Dernonzero}
  \Bigl(D^{(\alpha,\beta)}
  \, B^{\#}_{\alpha,\beta,\beta}\Bigr) ({-}1\,{,}\,{-}1) =
  \mfrac{\alpha!\,\beta!}{2^{\alpha+\beta}_{\nix}} \not= 0 \,,
\end{equation}
and we claim that for all other elements
$B^{\#}_{\widetilde\alpha,\widetilde\beta,\widetilde\gamma} \in
\widetilde\rI_k := \rI_k\setminus\{B^{\#}_{\alpha,\beta,\beta}\}$\,, we
have
\begin{equation}\label{eq:Derzero}
  \Bigl(D^{(n,m)}
  \, B^{\#}_{\widetilde\alpha,\widetilde\beta,\widetilde\gamma}\Bigr)
  ({-}1\,{,}\,{-}1) = 0
  \quad\text{for all}\quad(0\,,0)\leq(n,m)\leq(\alpha,\beta)\,.
\end{equation}
To this end, assume first that $\alpha<\beta$\,, then we have
\begin{equation*}
  \Bigl(D^{(n,m)}
  \, B^{\#}_{\beta,\alpha,\beta}\Bigr) ({-}1\,{,}\,{-}1) =0
  \quad\text{and}\quad
  \Bigl(D^{(n,m)}
  \, B^{\#}_{\beta,\beta,\alpha}\Bigr) ({-}1\,{,}\,{-}1) =0 \,,
\end{equation*}
since in~\eqref{eq:BoxspliDer}, we have
$\beta+\ell-(n-i)\geq\beta-\alpha>0$ and therefore we also have $\Bigl(
\mfrac{1{+}z_1}{2} \Bigr)^{\beta+\ell-(n-i)}=0$ for $z_1=-1$\,.
For arbitrary $\alpha\leq\beta$\,, consider
$(\alpha',\beta')\not=(\alpha,\beta)$\,.
  \par
In case $\alpha'<\alpha$ and thus $\beta'>\beta$\,, we have
\begin{equation}\label{eq:Derzeroabb}
  \Bigl(D^{(n,m)}
  \, B^{\#}_{\alpha',\beta',\beta'}\Bigr) ({-}1\,{,}\,{-}1) =0
\end{equation}
since in~\eqref{eq:BoxspliDer}, we have
$\beta'+\ell-(m-j)\geq\beta'-\beta>0$ and therefore we get $\Bigl(
\mfrac{1{+}z_2}{2} \Bigr)^{\beta'+\ell-(m-j)}=0$ for $z_2=-1$\,; and
also
\begin{equation}\label{eq:Derzerobabbba}
  \Bigl(D^{(n,m)}
  \, B^{\#}_{\beta',\alpha',\beta'}\Bigr) ({-}1\,{,}\,{-}1) =0
  \quad\text{and}\quad
  \Bigl(D^{(n,m)}
  \, B^{\#}_{\beta',\beta',\alpha'}\Bigr) ({-}1\,{,}\,{-}1) =0
\end{equation}
since in~\eqref{eq:BoxspliDer}, we have
$\beta'+\ell-(n-i)\geq\beta'-\alpha>0$ and therefore we also get $\Bigl(
\mfrac{1{+}z_1}{2} \Bigr)^{\beta'+\ell-(n-i)}=0$ for $z_1=-1$\,.
  \par
In case $\alpha<\alpha'\leq\beta'<\beta$\,, \eqref{eq:Derzeroabb} holds
since in~\eqref{eq:BoxspliDer}, we have
$\alpha'+\ell-(n-i)\geq\alpha'-\alpha>0$ and therefore $\Bigl(
\mfrac{1{+}z_1}{2} \Bigr)^{\alpha'+\ell-(n-i)}=0$ for $z_1=-1$\,; and
also \eqref{eq:Derzerobabbba} holds since
in~\eqref{eq:BoxspliDer}, we have
$\beta'+\ell-(n-i)\geq\beta'-\alpha>0$ and therefore $\Bigl(
\mfrac{1{+}z_1}{2} \Bigr)^{\beta'+\ell-(n-i)}=0$ for $z_1=-1$\,.
So \eqref{eq:Derzero} is shown.
  \par
But this implies that
\begin{equation*}
  B^{\#}_{\alpha,\beta,\beta} \notin
  {<}\,\widetilde\rI_k\,{>} \,,
\end{equation*}
since otherwise, we could write
\begin{equation*}
  B^{\#}_{\alpha,\beta,\beta} =
  \sum_{B^{\#}_{\widetilde\alpha,\widetilde\beta,\widetilde\gamma} \in
  \widetilde\rI_k}
  p_{\widetilde\alpha,\widetilde\beta,\widetilde\gamma} \,
  B^{\#}_{\widetilde\alpha,\widetilde\beta,\widetilde\gamma}
\end{equation*}
which, employing the Leibniz formula, yields
\begin{alignat*}{2}
  \Bigl(D^{(\alpha,\beta)} \, B^{\#}_{\alpha,\beta,\beta}\Bigr)
    &({-}1\,{,}\,{-}1) = \\
  = &\sum_{B^{\#}_{\widetilde\alpha,\widetilde\beta,\widetilde\gamma}
    \in \widetilde\rI_k} \sum_{n=0}^\alpha \mbinom{\alpha}{n}
    \sum_{m=0}^\beta \mbinom{\beta}{m} \,
    && \Bigl(D^{(\alpha-n,\beta-m)} \,
    p_{\widetilde\alpha,\widetilde\beta,\widetilde\gamma}\Bigr)
    ({-}1\,{,}\,{-}1) \\
  &&&\times \Bigl(D^{(n,m)} \,
    B^{\#}_{\widetilde\alpha,\widetilde\beta,\widetilde\gamma}\Bigr)
    ({-}1\,{,}\,{-}1) \,,
\end{alignat*}
and this contradicts \eqref{eq:Dernonzero} and \eqref{eq:Derzero}.
  \par
By symmetry in $z_1$ and $z_2$\,, it follows that also
\begin{equation*}
  B^{\#}_{\beta,\alpha,\beta} \notin
  {<}\,\rI_k\setminus\{B^{\#}_{\beta,\alpha,\beta}\}\,{>} \,.
\end{equation*}
  \par
It remains to show that
\begin{equation*}
  B^{\#}_{\beta,\beta,\alpha} \notin
  {<}\,\rI_k\setminus\{B^{\#}_{\beta,\beta,\alpha}\}\,{>} \,.
\end{equation*}
This can be achieved by employing a directional derivative and
considering mixed derivatives of the form
\begin{equation*}
  D^{(n,m,\ell)} = \Bigl(\mfrac{\partial}{\partial z_1}\Bigr)^n \;
  \Bigl(\mfrac{\partial}{\partial z_2}\Bigr)^m \;
  \Bigl(\mfrac{\partial}{\partial z_1} +
  \mfrac{\partial}{\partial z_2}\Bigr)^\ell =
  \sum_{j=0}^\ell \mbinom{\ell}{j} D^{(n+j,m+\ell-j)} \,.
\end{equation*}
Along the same lines as above, one shows that
\begin{equation*}
  \Bigl(D^{(\alpha,0,\beta)} \, B^{\#}_{\beta,\beta,\alpha}\Bigr)
  (1\,{,}\,{-}1) =
  (-1)^\alpha\mfrac{\alpha!\,\beta!}{2^{\alpha+\beta}_{\nix}} \not=0
  \,,
\end{equation*}
but that for all other elements
$B^{\#}_{\widetilde\alpha,\widetilde\beta,\widetilde\gamma} \in
\rI_k\setminus\{B^{\#}_{\beta,\beta,\alpha}\}$\,,
\begin{equation*}
  \Bigl(D^{(n,0,\ell)}
  \, B^{\#}_{\widetilde\alpha,\widetilde\beta,\widetilde\gamma}\Bigr)
  (1\,{,}\,{-}1) = 0
  \quad\text{for all }(0\,,0)\leq(n,\ell)\leq(\alpha,\beta)\,.
\end{equation*}
  \par
Alternatively, we may employ the coordinate transformation
\begin{equation*}
  \begin{pmatrix}w_1\\w_2\end{pmatrix} =
  \begin{pmatrix}z_1\,z_2\\1/z_2\end{pmatrix}
  \qquad\iff\qquad
   \begin{pmatrix}z_1\\z_2\end{pmatrix} =
  \begin{pmatrix}w_1\,w_2\\1/w_2\end{pmatrix}
\end{equation*}
which yields
\begin{equation*}
  B^{\#}_{\alpha,\beta,\gamma}(z_1,z_2) =
  \mfrac{1}{w_2^\beta}\,B^{\#}_{\gamma,\beta,\alpha}(w_1,w_2)
\end{equation*}
and thus allows us to reduce this case to the first one also.
\end{proof}

The characterization of convergent bivariate subdivision schemes
established in Theorem~C opens a way for their
systematic study. We only point out the following facts:
{\setlength{\leftmargini}{5mm}
   \begin{itemize}
\item
The set $\rI_k$ of generators for $\cI^k$ is symmetric in the sense
that it is invariant under an interchange of the two variables, and
that the indices of the generators can be permuted arbitrarily.
\item
For even $k$, the $k/2$-th power of (the symbol of) the Courant
hat function $B_{1,1,1}$ appears in $\rI_k$\,.

\smallskip
\item
Most interesting for us, however, is the fact that the smoothness
of these generators matches perfectly with the order of polynomial
reproduction. More generally, the box spline symbol
$B^{\#}_{\alpha,\beta,\gamma}$ is an element of
\begin{equation*}
  L_\infty^{(\kappa-1)} \subset C^{(\kappa-2)}
  \quad\text{for}\quad
  \kappa = \alpha + \beta + \gamma - \max\{\alpha,\beta,\gamma\} = k\,,
\end{equation*}
see \cite{BH2}. This means that for each three-directional box spline,
smoothness and
polynomial reproduction match in the same way as in the univariate
case. This property is not necessarily preserved under taking
combinations as described in Theorem~C, however, as
we shall see below for the four-directional box splines.

\smallskip
\item
It may also be of interest that the total degree of the box spline
with the symbol $B^\#_{\alpha,\beta,\gamma}$ is
\begin{equation*}
  \mu = \alpha+\beta+\gamma-2 \,,
\end{equation*}
which for the generators $B^\#_{\beta,\beta,\alpha}$\,,
$B^\#_{\beta,\alpha,\beta}$\,, and $B^\#_{\alpha,\beta,\beta}$ with
$\alpha\leq\beta$ yields
\begin{equation*}
  \mu = \alpha+2\,\beta-2 = 2\,k-\alpha-2 = 2\,\kappa-\alpha-2 \,.
\end{equation*}
\end{itemize}}

\subsection{Modification.}

As in Section~\ref{sec:zeroconideal}, we can replace the $q_\bTheta$ by
the $\widetilde{q}_\bTheta$\,, and all results still hold. More
precisely, this means replacing $B^{\#}_{\alpha,\beta,\gamma}$ by
\begin{equation*}
  \widetilde{B}^{\#}_{\alpha,\beta,\gamma}(z_1,z_2) =
  \Bigl(\mfrac{1+z_1}{2}\Bigr)^\alpha \,
  \Bigl(\mfrac{1+z_2}{2}\Bigr)^\beta \,
  \Bigl(\mfrac{z_1+z_2}{2}\Bigr)^\gamma \,.
\end{equation*}
For example, Lemma~\ref{lem:genred} holds as stated, and in its
proof, we only have to replace \eqref{eq:genred} by
\begin{equation*}
  \widetilde{B}_{1,0,0}^{\#}(z_1,z_2) +
  \widetilde{B}_{0,1,0}^{\#}(z_1,z_2) -
  \widetilde{B}_{0,0,1}^{\#}(z_1,z_2) = 1 \,.
\end{equation*}
Also, Theorem~\ref{thm:twoDJkgen} together with its proof,
Examples~\ref{exas:twoDJkgen} and Theorem~C are still
valid. Lemma~\ref{lem:BoxspliDer} now runs as follows.

\begin{lemma}\label{lem:modBoxspliDer}
The partial derivatives of the modified box spline symbol
$\widetilde{B}^{\#}_{\alpha,\beta,\gamma}$ are given by
\begin{align}
  \Bigl(D^{(n,m)} \,
    \widetilde{B}^{\#}_{\alpha,\beta,\gamma}\Bigr) &(z_1,z_2) =\notag\\
  = \sum_{\ell=0}^\gamma \mbinom{\gamma}{\ell}
    & \biggl( \sum_{i=0}^n \mfrac{n!}{2^n}
    \mbinom{\alpha}{n{-}i}
    \Bigl( \mfrac{1{+}z_1}{2} \Bigr)^{\alpha-(n-i)}
    \mbinom{\ell}{i}
    \Bigl( \mfrac{z_1}{2} \Bigr)^{\ell-i} \biggr)
    \label{eq:modBoxspliDer}\\
  &\times \biggl( \sum_{j=0}^m \mfrac{m!}{2^m}
    \mbinom{\beta}{m{-}j}
    \Bigl( \mfrac{1{+}z_2}{2} \Bigr)^{\beta-(m-j)}
    \mbinom{\gamma{-}\ell}{j}
    \Bigl( \mfrac{z_2}{2} \Bigr)^{\gamma-\ell-j} \biggr)\,. \notag
\end{align}
\end{lemma}

\begin{proof}
Writing
\begin{align*}
  \widetilde{B}^{\#}_{\alpha,\beta,\gamma} (z_1,z_2) &=
    \Bigl( \mfrac{1{+}z_1}{2} \Bigr)^\alpha
    \Bigl( \mfrac{1{+}z_2}{2} \Bigr)^\beta
    \Bigl( \mfrac{z_1}{2}+\mfrac{z_2}{2} \Bigr)^\gamma \\
 &= \sum_{\ell=0}^\gamma \mbinom{\gamma}{\ell}
   \Bigl( \mfrac{1{+}z_1}{2} \Bigr)^{\alpha}
   \Bigl( \mfrac{1{+}z_2}{2} \Bigr)^{\beta}
   \Bigl( \mfrac{z_1}{2} \Bigr)^{\ell}
   \Bigl( \mfrac{z_2}{2} \Bigr)^{\gamma-\ell}
\end{align*}
yields by the Leibniz formula
\begin{align*}
  \Bigl(D^{(n,m)}
    \, \widetilde{B}^{\#}_{\alpha,\beta,\gamma}\Bigr) & (z_1,z_2) = \\
  = \sum_{\ell=0}^\gamma \mbinom{\gamma}{\ell}
    & \mfrac{\partial^n}{\partial z_1^{\;n}}
    \biggl( \Bigl( \mfrac{1{+}z_1}{2} \Bigr)^{\alpha}
    \Bigl( \mfrac{z_1}{2} \Bigr)^{\ell} \biggr)
    \mfrac{\partial^m}{\partial z_2^{\;m}}
    \biggl( \Bigl( \mfrac{1{+}z_2}{2} \Bigr)^{\beta}
    \Bigl( \mfrac{z_2}{2} \Bigr)^{\gamma-\ell} \biggr) \\
  = \sum_{\ell=0}^\gamma \mbinom{\gamma}{\ell}
    & \biggl( \sum_{i=0}^n \mbinom{n}{i}
    \mfrac{\alpha!}{(\alpha{-}(n{-}i))!}
    \mfrac{(1{+}z_1)^{\alpha-(n-i)}}{2^{\alpha}_{\nix}}
    \mfrac{\ell!}{(\ell{-}i)!}
    \mfrac{z_1^{\ell-i}}{2^{\ell}_{\nix}} \biggr) \\
  &\times \biggl( \sum_{j=0}^m \mbinom{m}{j}
    \mfrac{\beta!}{(\beta{-}(m{-}j))!}
    \mfrac{(1{+}z_2)^{\beta-(m-j)}}{2^{\beta}_{\nix}}
    \mfrac{(\gamma{-}\ell)!}{(\gamma{-}\ell{-}j)!}
    \mfrac{z_2^{\gamma-\ell-j}}{2^{\gamma-\ell}_{\nix}} \biggr) \\
  = \sum_{\ell=0}^\gamma \mbinom{\gamma}{\ell}
    & \biggl( \sum_{i=0}^n \mfrac{n!}{2^n}
    \mbinom{\alpha}{n{-}i}
    \Bigl( \mfrac{1{+}z_1}{2} \Bigr)^{\alpha-(n-i)}
    \mbinom{\ell}{i}
    \Bigl( \mfrac{z_1}{2} \Bigr)^{\ell-i} \biggr) \\
  &\times \biggl( \sum_{j=0}^m \mfrac{m!}{2^m}
    \mbinom{\beta}{m{-}j}
    \Bigl( \mfrac{1{+}z_2}{2} \Bigr)^{\beta-(m-j)}
    \mbinom{\gamma{-}\ell}{j}
    \Bigl( \mfrac{z_2}{2} \Bigr)^{\gamma-\ell-j} \biggr)
\end{align*}
as claimed.
\end{proof}

This implies that Proposition~\ref{prp:maxJk} remains valid, and
in its proof, only the values of the listed nonzero mixed partial
derivatives need to be multiplied by $(-1)^\gamma$\,. Also,
Proposition~\ref{prp:mingen} still holds as stated: in its proof,
once more the signs of the nonzero derivative values change, and
in the last part, we need to consider the directional derivative
$\mfrac{\partial}{\partial z_1} - \mfrac{\partial}{\partial z_2}$
instead of $\mfrac{\partial}{\partial z_1} +
\mfrac{\partial}{\partial z_2}$\,; alternatively, we can
employ the coordinate transformation
\begin{equation*}
  \begin{pmatrix}w_1\\w_2\end{pmatrix} =
  \begin{pmatrix}z_1/z_2\\1/z_2\end{pmatrix}
  \qquad\iff\qquad
   \begin{pmatrix}z_1\\z_2\end{pmatrix} =
  \begin{pmatrix}w_1/w_2\\1/w_2\end{pmatrix}
\end{equation*}
which yields
\begin{equation*}
  \widetilde{B}^{\#}_{\alpha,\beta,\gamma}(z_1,z_2) =
  \mfrac{1}{w_2^{\alpha+\beta+\gamma}}\,
  \widetilde{B}^{\#}_{\gamma,\beta,\alpha}(w_1,w_2) \,.
\end{equation*}

\bigskip
Another approach is the use of the coordinate transformation
\begin{equation*}
  \begin{pmatrix}w_1\\w_2\end{pmatrix} =
  \begin{pmatrix}z_1\\1/z_2\end{pmatrix}
  \qquad\iff\qquad
  \begin{pmatrix}z_1\\z_2\end{pmatrix} =
  \begin{pmatrix}w_1\\1/w_2\end{pmatrix}
\end{equation*}
which was already mentioned at the beginning of
Section~\ref{sec:zeroconideal}. This yields
\begin{equation*}
  B^{\#}_{\alpha,\beta,\gamma}(z_1,z_2) =
  \mfrac{1}{w_2^{\beta+\gamma}}\,
  \widetilde{B}^{\#}_{\alpha,\beta,\gamma}(w_1,w_2)
\end{equation*}
which shows that all the results for the
$B^{\#}_{\alpha,\beta,\gamma}$ hold for the
$\widetilde{B}^{\#}_{\alpha,\beta,\gamma}$ also.

\subsection{Four-directional box splines}

The refinement mask symbol of a four-direc\-tional box spline has the
form
\begin{equation*}
  4\,B^{\#}_{\alpha,\beta,\gamma,\delta}(z_1,z_2) =
  4\,\Bigl(\mfrac{1+z_1}{2}\Bigr)^\alpha \,
  \Bigl(\mfrac{1+z_2}{2}\Bigr)^\beta \,
  \Bigl(\mfrac{1+z_1\,z_2}{2}\Bigr)^\gamma \,
  \Bigl(\mfrac{1+z_1/z_2}{2}\Bigr)^\delta \,.
\end{equation*}
This uses both alternatives, $\be_1+\be_2$ and $\be_1-\be_2$, for the third
direction discussed above. We can rewrite the above as
\begin{equation}\label{eq:fourdiralt}
  B^{\#}_{\alpha,\beta,\gamma,\delta}(z_1,z_2) =
  \mfrac{1}{2^\delta}\,\Bigl(1+\mfrac{z_1}{z_2}\Bigr)^\delta
  \,B^{\#}_{\alpha,\beta,\gamma}(z_1,z_2)\,,
\end{equation}
illustrating the well-known fact that any four-directional box
spline is indeed a special convex combination of the shifts of
some three-directional box spline. This convex combination uses
the normalized binomial weights $\lambda_{\delta,\ell} =
2^{-\delta}\,\binom{\delta}{\ell}$, $\ell=0,\dots,\delta$\,. The
representation in \eqref{eq:fourdiralt} is not optimal, however,
if one tries to determine the maximal $k$ such that
$B^\#_{\alpha,\beta,\gamma,\delta}\in\cI^k$\,. Instead, we use the
identity $z_1+z_2 = (1+z_1)\,(1+z_2) - (1+z_1\,z_2)$ which in
terms of the normalized box spline symbols reads as
\begin{equation*}
  B^\#_{0,0,0,1}(z_1,z_2) = \mfrac{1}{z_2}\,\mfrac{z_1+z_2}{2} =
  \mfrac{1}{z_2}\,\Bigl(2\,B^\#_{1,1,0}(z_1,z_2) -
  B^\#_{0,0,1}(z_1,z_2)\Bigr)\,.
\end{equation*}
This yields
\begin{align}\label{eq:fourdirid}
  B^\#_{\alpha,\beta,\gamma,\delta}(z_1,z_2)
    &= B^\#_{\alpha,\beta,\gamma}(z_1,z_2)\;\mfrac{1}{z_2^{\;\delta}}\,
    \Bigl( 2\,B^\#_{1,1,0}(z_1,z_2)-B^\#_{0,0,1}(z_1,z_2) \Bigr)^\delta
    \\
  &= \sum_{\ell=0}^\delta
    \mfrac{2^\ell\,(-1)^{\delta-\ell}}{z_2^{\;\delta}}\,
    \mbinom{\delta}{\ell}\,
    B^\#_{\alpha+\ell,\beta+\ell,\gamma+\delta-\ell}(z_1,z_2)\,.\notag
\end{align}

\begin{proposition}\label{prp:fourdirideal}
For any quadruple $(\alpha,\beta,\gamma,\delta)\in\NN_0^4$\,, the
maximal $k$ such that
$B^{\#}_{\alpha,\beta,\gamma,\delta}\in\cI^k$ is given by
\begin{equation}\label{eq:maxkfourdir}
  k = \alpha+\beta+\gamma+\delta -
  \max\{\alpha\,,\,\beta\,,\,\gamma{+}\delta\} \,.
\end{equation}
In other words, $B^{\#}_{\alpha,\beta,\gamma,\delta}\in\cI^k$ if and
only if $B^{\#}_{\alpha,\beta,\gamma+\delta}\in\cI^k$\,.
\end{proposition}

\begin{proof}
By Proposition~\ref{prp:maxJk}, we may conclude
from~\eqref{eq:fourdirid} that
$B^{\#}_{\alpha,\beta,\gamma,\delta}\in\cI^k$ for
\begin{align*}
  k &= \min_{\ell=0,\dots,\delta} \bigl\{
    \alpha+\beta+\gamma+\delta+\ell
    - \max \{ \alpha{+}\ell \,,\, \beta{+}\ell \,,\,
    \gamma{+}\delta{-}\ell \} \bigr\} \\
  &= \min_{\ell=0,\dots,\delta} \bigl\{ \min \{ \beta{+}\gamma{+}\delta
    \,,\, \alpha{+}\gamma{+}\delta \,,\, \alpha{+}\beta{+}2\,\ell \}
    \bigr\} \\
  &= \min \{ \beta{+}\gamma{+}\delta \,,\, \alpha{+}\gamma{+}\delta
    \,,\, \alpha{+}\beta \} \\
  &= \alpha+\beta+\gamma+\delta - \max\{\alpha \,,\, \beta \,,\,
    \gamma{+}\delta\} \,.
\end{align*}
Thus it remains to show that $B^{\#}_{\alpha,\beta,\gamma,\delta}
\notin \cI^{k+1}$ for this~$k$\,. Equivalently, we show that
\begin{equation*}
  \widetilde{B}^{\#}_{\alpha,\beta,\gamma,\delta}(z_1,z_2) =
  z_2^{\;\delta}\,B^{\#}_{\alpha,\beta,\gamma,\delta}(z_1,z_2) =
  \Bigl(\mfrac{1+z_1}{2}\Bigr)^\alpha \,
  \Bigl(\mfrac{1+z_2}{2}\Bigr)^\beta \,
  \Bigl(\mfrac{1+z_1\,z_2}{2}\Bigr)^\gamma \,
  \Bigl(\mfrac{z_1+z_2}{2}\Bigr)^\delta
\end{equation*}
satisfies $\widetilde{B}^{\#}_{\alpha,\beta,\gamma,\delta} \notin
\cI^{k+1}$ for this~$k$\,.

To this end, we proceed along the lines of
the proofs of Lemmata~\ref{lem:BoxspliDer} and~\ref{lem:modBoxspliDer}
and Proposition~\ref{prp:maxJk}. Writing
\begin{align*}
  \widetilde{B}^{\#}_{\alpha,\beta,\gamma,\delta} (z_1,z_2) &=
    \Bigl( \mfrac{1{+}z_1}{2} \Bigr)^\alpha
    \Bigl( \mfrac{1{+}z_2}{2} \Bigr)^\beta
    \Bigl( \mfrac{1+z_1}{2}\,\mfrac{1+z_2}{2} +
    \mfrac{z_1-1}{2}\,\mfrac{z_2-1}{2} \Bigr)^\gamma
    \Bigl( \mfrac{z_1}{2}+\mfrac{z_2}{2} \Bigr)^\delta \\
  &= \sum_{\ell=0}^\gamma \mbinom{\gamma}{\ell}
    \sum_{j=0}^\delta \mbinom{\delta}{j} \;
    \Bigl( \mfrac{1{+}z_1}{2} \Bigr)^{\alpha+\ell}
    \Bigl( \mfrac{z_1{-}1}{2} \Bigr)^{\gamma-\ell}
    \Bigl( \mfrac{z_1}{2} \Bigr)^{j} \\
  &\phantom{\;= \sum_{\ell=0}^\gamma \mbinom{\gamma}{\ell}
    \sum_{j=0}^\delta \mbinom{\delta}{j} \;}
    \times
    \Bigl( \mfrac{1{+}z_2}{2} \Bigr)^{\beta+\ell}
    \Bigl( \mfrac{z_2{-}1}{2} \Bigr)^{\gamma-\ell}
    \Bigl( \mfrac{z_2}{2} \Bigr)^{\delta-j}
\end{align*}
yields
\begin{align*}
  &\Bigl( D^{(n,m)} \,
  \widetilde{B}^{\#}_{\alpha,\beta,\gamma,\delta}\Bigr) (z_1,z_2) =  
    \sum_{\ell=0}^\gamma \mbinom{\gamma}{\ell} \\ & \hspace{6mm}\times
    \sum_{j=0}^\delta \mbinom{\delta}{j} \biggl( \sum_{\sum n_i=n} \mfrac{n!}{2^n}
    \mbinom{\alpha{+}\ell}{n_1}
    \Bigl( \mfrac{1{+}z_1}{2} \Bigr)^{\alpha+\ell-n_1}
    \mbinom{\gamma{-}\ell}{n_2}
    \Bigl( \mfrac{z_1{-}1}{2} \Bigr)^{\gamma-\ell-n_2}
    \mbinom{j}{n_3}
    \Bigl( \mfrac{z_1}{2} \Bigr)^{j-n_3} \biggr) \\
  &\hspace{6mm}\times \biggl( \sum_{\sum m_i=m} \mfrac{m!}{2^m}
    \mbinom{\beta{+}\ell}{m_1}
    \Bigl( \mfrac{1{+}z_2}{2} \Bigr)^{\beta+\ell-m_1}
    \mbinom{\gamma{-}\ell}{m_2}
    \Bigl( \mfrac{z_2{-}1}{2} \Bigr)^{\gamma-\ell-m_2}
    \mbinom{\delta{-}j}{m_3}
    \Bigl( \mfrac{z_2}{2} \Bigr)^{\delta-j-m_3} \biggr)
\end{align*}
and, therefore,
\begin{alignat*}{3}
  &\Bigl(D^{(\alpha,\beta)} \,
    \widetilde{B}^{\#}_{\alpha,\beta,\gamma,\delta}\Bigr)
    ({-}1\,{,}\,{-}1) &&=&
    (-1)^\delta\mfrac{\alpha!\,\beta!}{2^{\alpha+\beta}_{\nix}}
    &\not= 0 \,,\\
  &\Bigl(D^{(\alpha,\gamma+\delta)} \,
    \widetilde{B}^{\#}_{\alpha,\beta,\gamma,\delta}\Bigr)
    ({-}1\,{,}\,1) &&=& \;
    (-1)^\gamma
    \mfrac{\alpha!\,(\gamma{+}\delta)!}{2^{\alpha+\gamma+\delta}_{\nix}}
    &\not= 0 \,,\\
  \text{and}\quad
  &\Bigl(D^{(\gamma+\delta,\beta)} \,
    \widetilde{B}^{\#}_{\alpha,\beta,\gamma,\delta}\Bigr)
    (1\,{,}\,{-}1) &&=&
    (-1)^\gamma
    \mfrac{\beta!\,(\gamma{+}\delta)!}{2^{\beta+\gamma+\delta}_{\nix}}
    &\not= 0 \,.
\end{alignat*}
This shows that $\widetilde{B}^{\#}_{\alpha,\beta,\gamma,\delta} \notin
\cI^{k+1}$ for
\begin{equation*}
  k=\min\{\alpha{+}\beta\,,\,\alpha{+}\gamma{+}\delta\,,\,
  \beta{+}\gamma{+}\delta\} =
  \alpha+\beta+\gamma+\delta -
  \max\{\alpha\,,\,\beta\,,\,\gamma{+}\delta\}\,,
\end{equation*}
and this completes the proof.
\end{proof}

\medskip
\noindent {\bf Remark.} The proof of
Proposition~\ref{prp:fourdirideal} establishes the connection
between the order $k$ in \eqref{eq:maxkfourdir} and the smoothness
of the $(\alpha,\beta,\gamma,\delta)$-box spline, which can be
determined using the results in~\cite{BH1}. The function is an
element of $L_\infty^{(\kappa-1)} \subset C^{(\kappa-2)}$ for
\begin{equation*}
  \kappa = \min \{\alpha{+}\beta{+}\gamma \,,\, \alpha{+}\beta{+}\delta
  \,,\, \alpha{+}\gamma{+}\delta \,,\, \beta{+}\gamma{+}\delta \} =
  \alpha+\beta+\gamma+\delta -
  \max\{\alpha\,,\,\beta\,,\,\gamma\,,\,\delta\}\,,
\end{equation*}
while the order $k$ (degree $k-1$) of polynomial reproduction is
given by
\begin{equation*}
  k=\min\{\alpha{+}\beta\,,\,\alpha{+}\gamma{+}\delta\,,\,
  \beta{+}\gamma{+}\delta\} =
  \alpha+\beta+\gamma+\delta -
  \max\{\alpha\,,\,\beta\,,\,\gamma{+}\delta\}\,.
\end{equation*}

\begin{corollary}
For the four-directional $(\alpha,\beta,\gamma,\delta)$-box
spline, we have $\kappa\geq k$\,, and
\begin{equation*}
  \kappa > k \quad\iff\quad \gamma+\delta > \max\{\alpha,\beta\}
  \quad\text{and}\quad \min\{\gamma,\delta\} > 0\,.
\end{equation*}
\end{corollary}
\begin{proof}
By the above remark,
\begin{equation*}
  \kappa-k = \max\{\alpha\,,\,\beta\,,\,\gamma{+}\delta\} -
  \max\{\alpha\,,\,\beta\,,\,\gamma\,,\,\delta\} \ge 0 \;.
\end{equation*}
In the case
$\max\{\alpha,\beta,\gamma{+}\delta\} = \max\{\alpha,\beta\}$\,, we
have $\max\{\gamma,\delta\} \leq \gamma{+}\delta \leq
\max\{\alpha,\beta\}$\,, \textit{i.\,e.}, also
$\max\{\alpha,\beta,\gamma,\delta\} = \max\{\alpha,\beta\}$ and
therefore $\kappa=k$\,. Otherwise, we have $\gamma{+}\delta >
\max\{\alpha,\beta\}$\,, and then
\begin{equation*}
  \kappa-k = \gamma{+}\delta - \max\{\alpha,\beta,\gamma,\delta\} =
  \min\{\gamma{+}\delta{-}\alpha\,{,}\,\gamma{+}\delta{-}\beta\,{,}\,
  \gamma\,{,}\,\delta\}
\end{equation*}
where the first two elements are positive, so in this case,
\begin{equation*}
  \kappa-k = 0 \quad\iff\quad \min\{\gamma,\delta\} = 0\,.
\end{equation*}
All in all, this yields that
\begin{equation*}
  \kappa-k = 0 \qquad\iff\qquad \gamma{+}\delta\leq\max\{\alpha,\beta\}
  \quad\text{or}\quad \min\{\gamma,\delta\} = 0\,,
\end{equation*}
which proves the claim.
\end{proof}

The condition $\min\{\gamma,\delta\} = 0$ is worth a closer look. For
$\delta=0$\,, we are in the standard three-directional case, while for
$\gamma=0$\,, we have
\begin{equation*}
  B^{\#}_{\alpha,\beta,0,\delta}(z_1,z_2) =
  z_2^{\,\beta}\,B^{\#}_{\alpha,\beta,\delta}(z_1,1/z_2) =
  z_2^{-\delta}\,\widetilde{B}^{\#}_{\alpha,\beta,\delta}(z_1,z_2)\,,
\end{equation*}
so these are three-directional splines on the reflected grid.

\begin{examples}
We list a few standard examples together with the decompositions of
their mask symbols according to~\eqref{eq:fourdirid}
and to~\eqref{eq:bivgen}.

The $(1,1,1,1)$-box spline, known as the Zwart-Powell element,
has the mask symbol
\begin{align*}
  4\,B^{\#}_{1,1,1,1}(z_1,z_2) &= \mfrac{4}{z_2}\,
    \Bigl( 2\,B^{\#}_{2,2,1}(z_1,z_2) - B^{\#}_{1,1,2}(z_1,z_2) \Bigr) \\
  &= 1\cdot\mfrac{1+z_1/z_2}{2}\cdot4\,B^{\#}_{1,1,1}(z_1,z_2)\,.
\end{align*}
We find $k = 4-2=2$\,, so the associated subdivision scheme reproduces
polynomials of total degree up to one; but $\kappa=4-1=3$\,, so the
function is in $L_\infty^{(2)} \subset C^{(1)}$\,.

\medskip
The $(2,2,1,1)$-box spline has the mask symbol
\begin{align*}
  4\,B^{\#}_{2,2,1,1}(z_1,z_2) &= \mfrac{4}{z_2}\,
    \Bigl( 2\,B^{\#}_{3,3,1}(z_1,z_2) - B^{\#}_{2,2,2}(z_1,z_2) \Bigr) \\
  &= 2\cdot\mfrac{1}{z_2}\cdot4\,B^{\#}_{3,3,1}(z_1,z_2) +
    (-1)\cdot\mfrac{1}{z_2}\cdot4\,B^{\#}_{2,2,2}(z_1,z_2) \,.
\end{align*}
Here, we have $k = 6-2=4$\,, telling us that polynomials of degree up
to $3$ are reproduced, and also  $\kappa=6-2=4$\,, \textit{i.\,e.}, the
function is an element of $L_\infty^{(3)} \subset C^{(2)}$\,.

\medskip
More interesting are the higher order four-directional splines. For
example, the $(4,4,1,1)$-box spline has order of polynomial
reproduction $k=10-4=6$ and $\kappa=10-4=6$\,. Its mask symbol can be
represented as
\begin{align*}
  4\,B^{\#}_{4,4,1,1}(z_1,z_2) &= \mfrac{4}{z_2}\,
    \Bigl( 2\,B^{\#}_{5,5,1} - B^{\#}_{4,4,2}(z_1,z_2) \Bigr) \\
  &= 2\cdot\mfrac{1}{z_2}\cdot4\,B^{\#}_{5,5,1}(z_1,z_2) +
    (-1)\cdot\mfrac{1}{z_2}\cdot4\,B^{\#}_{4,4,2}(z_1,z_2) \,.
\end{align*}
\end{examples}

\medskip
\begin{remark}
The vector case, {\it i.\,e.}, the case when the matrix mask $\rA \in
\ell_0^{s \times s}(\ZZ^d)$ is a finitely supported sequence of $s
\times s$-matrices indexed by $\ZZ^d$, is more intricate. The
formulation of the zero conditions~\eqref{eq:Zkcond} for the mask
symbol of multivariate vector subdivision schemes depends greatly on
the so-called rank of the scheme, see \cite{ChaConSau05, CotroneiSauer,
JJL}. Such a formulation, see, {\it e.\,g.}, \cite{H03}, does not allow
us to read off the properties of the entries of the matrix Laurent
polynomial $\rA(\bz)$ directly. It is possible, though, to use a slight
modification of the transformation in \cite{H03} to obtain the matrix
sequences $\rT, \rT^{\textit{inv}} \in \ell_0^{s \times s}(\ZZ^d)$ such
that
\begin{equation*}
  \widetilde{\rA}(\bz)= \rT^{\textit{inv}}(\bz^2) \cdot \rA (\bz) \cdot
  \rT(\bz)
\end{equation*}
satisfies the zero conditions of a form that makes the structure of
some of the entries of $\rA(\bz)$ more evident. However, this is a
topic for further investigations.
\end{remark}

\section{Examples} 

\noindent In this section we illustrate the result of
Theorem~C with some examples. We would like to
emphasize that this result does not only simplify the study of the
properties of subdivision schemes, but also yields a way for
enhancing certain properties of existing schemes by combining them
appropriately.

In the following, the set $\rI_k$ of the box spline symbols is the
set of generators for $\cI^k$ as in Theorem~\ref{thm:twoDJkgen}.

In the masks displayed below the boldface entry at bottom-left
position refers to the index $(0,0)$. This assumption is not
really important, but as stated above already, we prefer to shift
masks so that they are supported in the positive quadrant and have
polynomial symbols.

\subsection{A bivariate interpolatory scheme}

Interpolatory schemes are characterized by the fact that one of
the submasks is a $\delta$ sequence, or equivalently, one of the
subsymbols is identically one. Theorem~C allows us to present a
systematic way for creating interpolatory schemes from our lists
of generators by equating the coefficients of their affine
combinations and normalizing them appropriately.

To provide just one such example, consider the interpolatory
scheme studied in \cite[Example 2]{Han2003}, a bivariate version of the
univariate four-point interpolation scheme given in
\cite{DynGreLev}. Its mask is
\begin{equation*}
  \ra = \frac{1}{32} \left(\begin{matrix}
          0 &  0 & -1 & -2 & -1 &  0 &  0 \\
          0 &  0 &  0 &  0 &  0 &  0 &  0 \\
         -1 &  0 & 10 & 18 & 10 &  0 & -1 \\
         -2 &  0 & 18 & 32 & 18 &  0 & -2 \\
         -1 &  0 & 10 & 18 & 10 &  0 & -1 \\
          0 &  0 &  0 &  0 &  0 &  0 &  0 \\
    {\bf 0} &  0 & -1 & -2 & -1 &  0 &  0 \\
  \end{matrix} \right)\;.
\end{equation*}
The scheme reproduces polynomials up to degree $k-1=3$, whence
$\ra(\bz)\in\cI^4$, and a representation of $\ra(\bz)$ in terms of
three-directional box splines from the list $\rI_4$ is given by
\begin{align*}
  &z_1^3 z_2^3\;\ra(z_1,z_2) \\
  &= \; -2^4\,B^\#_{4,4,0}(z_1,z_2)
        - 2\,(z_1^2\,{+}\,z_2^2)\,B^\#_{2,2,2}(z_1,z_2)
        + 2^3\,(1\,{+}\,z_1\,{+}\,z_2)\, B^\#_{3,3,1}(z_1,z_2) \\
  &= \; 4\,\Big\{ {-}4\,B^\#_{4,4,0}(z_1,z_2)
                   - \mfrac{z_1^2+z_2^2}{2}\,B^\#_{2,2,2}(z_1,z_2)
                   + 6\,\mfrac{1+z_1+z_2}{3}\,B^\#_{3,3,1}(z_1,z_2)
           \Big\}\,.
\end{align*}
From the second line, the weights $\lambda$ are recognized as
$-4$, $-1$, and $6$, and the normalized $\sigma$-symbols are
\begin{equation*}
  1\;, \quad \frac{z_1^2+z_2^2}{2}\;, \quad\text{and}\quad
  \frac{1+z_1+z_2}{3}\;,
\end{equation*}
respectively.

\subsection{The butterfly scheme}\label{sec:butterfly}

The butterfly scheme has been studied in \cite{DynLevMic} and
\cite[ Example 5]{Han2003}. Its mask is given by
\begin{equation*}
  \ra = \frac{1}{16} \left(\begin{matrix}
          0 &  0 &  0 &  0 & -1 & -1 &  0  \\
          0 &  0 & -1 &  0 &  2 &  0 & -1  \\
          0 & -1 &  2 &  8 &  8 &  2 & -1  \\
          0 &  0 &  8 & 16 &  8 &  0 &  0  \\
         -1 &  2 &  8 &  8 &  2 & -1 &  0  \\
         -1 &  0 &  2 &  0 & -1 &  0 &  0  \\
    {\bf 0} & -1 & -1 &  0 &  0 &  0 &  0  \\
  \end{matrix} \right)\;.
\end{equation*}
It is  an interpolating scheme, and reproduces polynomials of
degree $k-1=3$. The representation of the mask symbol in terms of
three-directional box spline symbols from the list $\rI_4$ is
given by
\begin{align*}
  z_1^3 z_2^3&\,\ra(z_1,z_2) \\
  &= 4\, \Big\{ 26\,\mfrac{7+6\,z_1 z_2}{13}\,B^\#_{3,3,1}(z_1,z_2)
    - 2\,z_2 B^\#_{3,1,3}(z_1,z_2)
    - 2\,z_1 B^\#_{1,3,3}(z_1,z_2) \\
  & \hspace{60mm} -21 \, \mfrac{1+z_1+z_2}{3}\,B^\#_{2,2,2}(z_1,z_2)
    \Big\} \\
  &= 4\, \Big\{ 7\,z_1 z_2 B^\#_{2,2,2}(z_1,z_2) \\
  &\hspace{20mm}
    \,{-}\, 2\,z_1 B^\#_{1,3,3}(z_1,z_2)
    \,{-}\, 2\,z_2 B^\#_{3,1,3}(z_1,z_2)
    \,{-}\, 2\,z_1 z_2 B^\#_{3,3,1}(z_1,z_2) \Big\}\,.
\end{align*}
We see that the generators are all multiples of
$B^\#_{1,1,1}(z_1,z_2)$. This tells us that the symbol can be
factorized as
\begin{equation*}
  \ra (z_1,z_2) = B^\#_{1,1,1}(z_1,z_2)\; \rb(z_1,z_2)\;,
\end{equation*}
a fact noticed in \cite{DynLevMic}. We would like to emphasize the
following properties of the butterfly scheme. Firstly, a simple
computation yields that the symbol $\rb(\bz)$ does not define a
convergent subdivision scheme, although each of the summands in
$\rb(\bz)$ by itself does correspond to a convergent scheme. Secondly,
butterfly is an interpolatory subdivision scheme, but none of the
summands in the affine combination above possess this property.


\subsection{A convergent scheme}
The symbols presented in the above examples all possess a property
that is very important for their regularity analysis: they are
multiples of one specific box spline symbol of type
$B^\#_{\alpha,\beta,\gamma}$\,. The regularity analysis of such
schemes is given in \cite[ Section~4.3]{DynLevin02}. The type of
factorization used there, however, is a very special situation
which does not generally hold for convergent schemes.

A very simple example that comes to mind is the symbol given by
\begin{align*}
  \ra(z_1,z_2)&= 4\; \left\{\frac{1}{2}\, B^\#_{1,1,0}(z_1,z_2) +
    \frac{1}{2}\, B^\#_{0,1,2}(z_1,z_2)\right\} =
    4\;\frac{1+z_2}{2} \;\rc(z_1,z_2) \\
  \text{with}\qquad
    \rc(z_1,z_2) &= \frac{1}{2}\;\frac{1+z_1}{2} + \frac{1}{2}\;
    \left( \frac{1+z_1z_2}{2}\right)^2 \;.
\end{align*}
By Theorem~\ref{thm:twoDJkgen}, the symbol $\ra(\bz)$ is in
$\cI$, but none of the generators from the list $\rI_1$ divides
the symbol.

In order to check the convergence of this scheme, we study the properties of
the so-called difference scheme $\cS_\rB$, see  \cite{CDM91, Sauer},
with the matrix mask symbol $\rB(\bz)$ satisfying
\begin{equation*}
  \ra(\bz)\;\left(\begin{matrix} 1-z_1 \\ 1-z_2 \end{matrix} \right)^T =
  \left(\begin{matrix} 1-z_1^2 \\ 1-z_2^2
  \end{matrix} \right)^T \;
 \rB(\bz)\;.
\end{equation*}
One possible
such $\rB(\bz)$ is given by
\begin{alignat*}{2}
  && \rB(z_1,z_2)    &= \begin{pmatrix}
                          \rb_{11}(z_1,z_2) & \rb_{12}(z_1,z_2) \\
                          \rb_{21}(z_1,z_2) & \rb_{22}(z_1,z_2)
                        \end{pmatrix} \\[2mm]
  &\text{with} \qquad&
     \rb_{11}(z_1,z_2) &=
     \frac{1}{4}\; (z_1z_2^3-z_2^3+z_1z_2^2+z_2^2+4z_2+2 )\;,  \\
  && \rb_{12}(z_1,z_2) &= 0 \;, \\
  && \rb_{21}(z_1,z_2) &= \frac{1}{4}\; (z_1z_2-z_1-z_2+1 ) \;, \\
  &\text{and}&
     \rb_{22}(z_1,z_2) &= \frac{1}{4}\; (z_1^2z_2^2+2z_1z_2+2z_1+3 )\;.
\end{alignat*}
To check the convergence, we have to verify that the vector subdivision scheme $\cS_\rB$
converges to zero, see \cite{ChaConSau05,D92}. The   symbolic
calculations yield $\|\cS_\rB^5\|_\infty <1$ and, thus, that $\cS_\ra$ is
$C$-convergent.

It is worth noting that in this example the two building blocks,
with symbols $4\;B^\#_{1,1,0}(z_1,z_2)$ and
$4\;B^\#_{0,1,2}(z_1,z_2)$, are not the symbols of $C$-convergent
subdivision schemes, while the combination yields
$C$-convergence.

We also refer to the constructions
in \cite{CGPS07}, where the convex combination of a four-directional,
zero order box spline and a $C^1$-quadratic box spline are used to
obtain the so-called GP pseudo-quadratic box spline. This example shows
enhancement with respect to linear independence of the translates, at
the expense of reduced joint smoothness.

\bibliographystyle{amsplain}


\end{document}